\documentclass[12pt]{amsart}
\usepackage{graphicx,color}
\usepackage{a4wide}
\usepackage{amsthm}
\usepackage{amsopn}
\usepackage{amssymb}
\usepackage{amsmath}
\usepackage{latexsym}
\usepackage[utf8]{inputenc}
\usepackage[T1]{fontenc}
\usepackage{amsfonts}
\usepackage[all,cmtip]{xy}
\usepackage{mathrsfs}
\usepackage{tikz}
\usepackage{tikz-cd}
\usetikzlibrary{matrix, calc, arrows,backgrounds}
\usepackage[colorlinks=true, linkcolor=blue, citecolor=blue, urlcolor=blue]{hyperref}

\tikzset{node distance=2cm,auto}

\newtheorem{theorem}{Theorem}[section]
\newtheorem{lemma}[theorem]{Lemma}
\newtheorem{corollary}[theorem]{Corollary}
\newtheorem{proposition}[theorem]{Proposition}
\newtheorem{conjecture}{Conjecture}

\newenvironment{namedtheorem}[1]{\hfill\\ \noindent{\bf#1} \ \it}{\hfill\\}
\newtheorem{mainthm}{Theorem}

\newtheorem{maincor}[mainthm]{Corollary}

\theoremstyle{definition}
\newtheorem{question}[theorem]{Question}
\newtheorem{definition}[theorem]{Definition}

\newtheorem{assumption}{Assumption}
\newtheorem{example}[theorem]{Example}

\newtheorem{remark}[theorem]{Remark}

\usepackage{amsfonts}
\usepackage{color}

\newcommand{\bpf}{\noindent{\bf Proof}\hspace{7pt}}
\newcommand{\epf}{\qed}
\newcommand{\ben}{\begin{enumerate}}
\newcommand{\een}{\end{enumerate}}
\newcommand{\ble}{\begin{lemma}}
\newcommand{\ele}{\end{lemma}}
\newcommand{\bth}{\begin{theorem}}
\renewcommand{\eth}{\end{theorem}}
\newcommand{\bmth}{\begin{mainthm}}
\newcommand{\emth}{\end{mainthm}}
\newcommand{\bmco}{\begin{maincor}}
\newcommand{\emco}{\end{maincor}}

\newcommand{\bpr}{\begin{proposition}}
\newcommand{\epr}{\end{proposition}}
\newcommand{\bco}{\begin{corollary}}
\newcommand{\eco}{\end{corollary}}
\newcommand{\bcon}{\begin{conjecture}}
\newcommand{\econ}{\end{conjecture}}
\newcommand{\bqu}{\begin{question}}
\newcommand{\equ}{\end{question}}
\newcommand{\bde}{\begin{definition}}
\newcommand{\ede}{\end{definition}}
\newcommand{\bas}{\begin{assumption}}
\newcommand{\eas}{\end{assumption}}
\newcommand{\bre}{\begin{remark}}
\newcommand{\ere}{\end{remark}}
\newcommand{\bex}{\begin{example}}
\newcommand{\eex}{\end{example}}
\newcommand{\barr}{\begin{array}}
\newcommand{\earr}{\end{array}}
\newcommand{\btab}{\begin{tabular}}
\newcommand{\etab}{\end{tabular}}
\newcommand{\beq}{\begin{equation}}
\newcommand{\eeq}{\end{equation}}
\newcommand{\bea}{\begin{eqnarray*}}
\newcommand{\eea}{\end{eqnarray*}}
\newcommand{\bce}{\begin{center}}
\newcommand{\ece}{\end{center}}
\newcommand{\bpi}{\begin{picture}}
\newcommand{\epi}{\end{picture}}
\newcommand{\bfi}{\begin{figure} \begin{center}}
\newcommand{\efi}{\end{center} \end{figure}}

\newcommand{\bsl}{\begin{slide}{}}
\newcommand{\esl}{\end{slide}}

\newcommand{\pf}{\noindent{\bf Proof}\hspace{7pt}}

\newcommand{\hso}[1]{\hspace{-1pt}}

\newcommand{\into}{\hookrightarrow}

\newcommand{\sbe}{\subseteq}
\newcommand{\sbne}{\subsetneq}

\def\<{\langle}
\def\>{\rangle}

\newcommand{\De}{\Delta}

\newcommand{\La}{\Lambda}

\newcommand{\cB}{{\mathcal B}}

\newcommand{\cP}{{\mathcal P}}

\newcommand{\hI}{{I(\hGm)}}
\newcommand{\hJ}{{\widehat{J}}}

\newcommand{\Aut}{\mathop{\rm Aut}\nolimits}

\renewcommand{\bar}{\overline}

\newcommand{\End}{\mathop{\rm End}\nolimits}

\newcommand{\id}{\mathop{\rm id}\nolimits}
\newcommand{\im}{\mathop{\rm im}\nolimits}

\def\flexbox#1{\mathchoice{\mbox{#1}}{\mbox{#1}}{\mbox{\scriptsize #1}}%
{\mbox{\tiny #1}}}
\def\SL{\mathop{\flexbox{\rm SL}}\nolimits}
\def\PSL{\mathop{\flexbox{\rm PSL}}\nolimits}
\newcommand{\POmega}{\mathop{\flexbox{\rm P}\Omega}\nolimits}

\newcommand{\Char}{\mathop{\flexbox{\rm Char}}}

\newcommand{\FF}{{\mathbb F}}

\newcommand{\NN}{{\mathbb N}}

\newcommand{\ZZ}{{\mathbb Z}}

\newcommand{\comop}{{\mathscr K}}
\newcommand\hGm{\Lambda}
\newcommand{\hE}{{E(\hGm)}}

\newcommand{\hPi}{{\Pi}}

\newcommand{\KM}{G}
\newcommand{\compLth}{{\compL^\theta}}

\newcommand{\Pig}{\pi(\Gamma,{0})}

\newcommand{\tH}{{\tilde{H}}}

\newcommand\Gm{\Gamma}

\newcommand{\opp}{\mathbin{\rm opp}}

\newcommand{\fk}{{\mathsf k}}

\newcommand{\dfn}{\em}

\newcommand{\after}{\mathbin{ \circ }}

\newcommand{\Stab}{\mathop{\rm Stab}}

\newcommand{\normal}{\lhd}
\DeclareMathOperator{\proj}{proj}

\renewcommand{\hat}{\widehat}

\newcommand{\vep}{\varepsilon}

\renewcommand{\qed}{\hfill $\square$}

\newcounter{romanlistctr}
{\end{list}}%

 \makeatletter
 \@addtoreset{equation}{section}
 \makeatother

 \makeatletter
 \def\section{\@startsection {section}{1}{\z@}{-1.5ex plus -.5ex
 minus -.2ex}{1ex plus .2ex}{\large\bf}}

 \makeatletter
 \def\subsection{\@startsection {subsection}{1}{\z@}{-1.5ex plus -.5ex
 minus -.2ex}{1ex plus .2ex}{\bf}}

\newcommand{\amgrpA}{{\mathbf{A}}}
\newcommand{\amgrpD}{{\bf{D}}}
\newcommand{\amgrpG}{{\mathbf{G}}}
\newcommand{\amgrpL}{{\mathbf{L}}}
\newcommand{\amgrpU}{{\mathbf{U}}}

\newcommand{\ama}{{\mathbf a}}
\newcommand{\amg}{{\mathbf g}}

\newcommand{\aml}{{\mathbf l}}
\newcommand{\amz}{{\mathbf z}}

\newcommand{\ucompmapf}{{\tilde{\phi}}}

\newcommand{\ammapf}{\phi}

\newcommand{\ammapm}{\mu}
\newcommand{\ammapp}{\pi}
\newcommand{\uammapp}{\tilde{\pi}}

\newcommand{\compA}{{{A}}}
\newcommand{\compD}{{{D}}}
\newcommand{\compG}{{{G}}}
\newcommand{\compL}{{{L}}}
\newcommand{\compN}{{{N}}}

\newcommand{\compB}{{{B}}}
\newcommand{\compU}{{{U}}}
\newcommand{\compa}{\alpha}

\newcommand{\compl}{\lambda}

\newcommand{\compg}{\gamma}

\newcommand{\ucompA}{{\tilde{A}}}
\newcommand{\ucompB}{{\tilde{B}}}
\newcommand{\ucompD}{{\tilde{D}}}
\newcommand{\ucompG}{{\tilde{G}}}
\newcommand{\ucompH}{{\tilde{H}}}
\newcommand{\ucompK}{{\tilde{K}}}
\newcommand{\ucompL}{{\tilde{L}}}

\newcommand{\ucompU}{{\tilde{U}}}
\newcommand{\ucompZ}{{\tilde{Z}}}

\newcommand{\ucompz}{{\tilde{z}}}

\newcommand{\ucompa}{{\tilde{\alpha}}}

\newcommand{\ucompg}{{\tilde{\gamma}}}
\newcommand{\ucompl}{{\tilde{\lambda}}}

\newcommand{\ocompD}{{{D}}} %oriented completion
\newcommand{\ocompN}{{{N}}} %oriented completion

\newcommand{\umap}{\pi}

\newcommand{\kIz}{{\mathsf{z}}}
\newcommand{\kIF}{{\mathsf{F}}}
\newcommand{\kIK}{{\mathsf{K}}}
\newcommand{\kIM}{{\mathsf{M}}}
\newcommand{\imnu}{\kIM}

\newcommand{\kIZ}{{\mathsf{Z}}}

\newcommand{\rsc}{\mathrm{sc}}

\def\<{\langle}
\def\>{\rangle}

\newcommand{\tamB}{\tilde{\mathscr{B}}}
\newcommand{\amA}{{\mathscr{A}}}

\newcommand{\amG}{{\mathscr{G}}}

\newcommand{\amL}{{\mathscr{L}}}

\newcommand{\amB}{{\mathscr{B}}}
\DeclareMathOperator{\Spin}{Spin}

\newcommand{\classmap}{\omega}

\newcommand{\nom}{}
\newcommand{\idx}{}

\newcommand{\ead}{\email}
\newenvironment{keyword}{Keywords: \keywords}{}

\newcommand{\MSC}[1]{MSC #1:\subjclass}
\newcommand{\sep}{\hspace{1ex}}
\begin{document}
%\begin{frontmatter}
\title{\Large Curtis-Tits Groups of simply-laced type}

\author{Rieuwert J. Blok}
\ead{blokr@member.ams.org}
\address{Department of Mathematics and Statistics\\
Bowling Green State University\\
Bowling Green, OH 43403\\
U.S.A.}
\author{Corneliu G. Hoffman}
\ead{C.G.Hoffman@bham.ac.uk}
\address{University of Birmingham\\
Edgbaston, B15 2TT\\
U.K.}

\begin{abstract}
The classification of Curtis-Tits amalgams with {connected}, triangle free, simply-laced diagram over a field of size at least $4$ was completed in~\cite{BloHof2014b}. Orientable amalgams are those arising from applying the Curtis-Tits theorem to groups of Kac-Moody type, and indeed, their universal completions are central extensions of those groups of Kac-Moody type. 
The paper~\cite{BloHof2014a} exhibits concrete (matrix) groups as completions for all Curtis-Tits amalgams with diagram $\widetilde{A}_{n-1}$. For non-orientable amalgams these groups are symmetry groups of certain unitary forms over a ring of skew Laurent polynomials.
In the present paper we generalize this to all amalgams arising from the classification above and, under some additional conditions, exhibit their universal completions as central extensions of twisted groups of Kac-Moody type.
 \end{abstract}
\maketitle

%\end{frontmatter}

\section{Introduction}
 A celebrated theorem of Curtis and Tits~\cite{Cur1965a,Ti1974} (later extended by  Timmesfeld (see \cite{Tim1998,Tim03,Tim04,Tim06} for spherical groups), by P. Abramenko and  B.~M\"uhlherr~\cite{AbrMuh97} and Caprace~\cite{Cap2007} to 2-spherical groups { of Kac-Moody type}) on groups with finite BN-pair states that { these groups} are central quotients of the universal completion of the concrete amalgam of the Levi components of the parabolic subgroups with respect to a given (twin-) BN-pair.
Following Tits~\cite{Ti1992} a \nom{{\em group of Kac-Moody type}}{} is by definition a group with RGD system such that a central quotient is the subgroup of $\Aut(\De)$ generated by the root groups of an apartment in a Moufang twin-building \nom{$\De$}{}. This central quotient will be called the associated {\em adjoint} group of Kac-Moody type.

\subsection{Classification of Curtis-Tits amalgams}
In~\cite{BloHof2014b}, for a given {connected} diagram, we consider all possible rank-$2$ amalgams which are locally isomorphic to an amalgam arising from the Curtis-Tits theorem { in the simply-laced case}. 
Then,  extending the techniques from~\cite{BloHof2014a}, we classify Curtis-Tits amalgams as defined in Definition~\ref{dfn:CT structure} (called Curtis-Tits structures in loc.~cit.) with \nom{property (D)}{} (see Definition~\ref{dfn:property D}) and connected simply-laced diagram. We show that any Curtis-Tits amalgam with connected simply-laced diagram without triangles that has a non-trivial universal completion must have property (D). 
Although many of the results in the present and preceding papers can be proved under the more general assumption of property (D), some rely on $3$-sphericity, for instance the simple connectedness of $\De^\theta$ in Section~\ref{section:De theta} and Lemma~\ref{lem:vertex and edge group images}.

\bas
Throughout the paper we will make the following assumptions:
\begin{enumerate}
\item \nom{$\fk$}{} is a field with at least $4$ elements.
\item  \nom{$\Gamma=(I,E)$}{} is a connected simply-laced {Dynkin} diagram of rank $n=|I|$ without triangles.
\item  All Curtis-Tits amalgams have property (D).
\end{enumerate}
\eas

In this context, we quote Theorem 1 from~\cite{BloHof2014b}:

\begin{namedtheorem}{Classification Theorem}\label{thm:CT amalgams} Let $\Gamma$ be a connected simply laced {Dynkin} diagram without triangles and 
{$\fk$}{} a field with at least 4 elements.
There is a natural bijection between isomorphism classes of universal Curtis-Tits amalgams with property (D) over the field $\fk$ on the graph $\Gamma$  and group homomorphisms \nom{$\classmap\colon\Pig\to\langle\tau\rangle\times\Aut(\fk)$}{}.
 \end{namedtheorem}
 
 We'd like to point out that the classification theorem generalizes similar classification results of  sound Moufang foundations in~\cite{Mu1999,Ti1992}. 
 
Here, \nom{$\Pig$}{} denotes the (first) fundamental group of the graph $\Gamma$ with base point $0\in I$, and \nom{$\langle\tau\rangle\times\Aut(\fk)$}{} arises as a group of automorphisms of $\SL_2(\fk)$ stabilizing a given torus (see Subsection~\ref{subsec:universal CT amalgams}).
It is clear from their definition that any Curtis-Tits amalgam is a homomorphic image of a universal Curtis-Tits amalgam. Note that the standard Curtis-Tits amalgam of type $\Gamma_\classmap(\fk)$ in Definition~\ref{dfn:standard CT structure} is universal. We shall denote this amalgam by \nom{$\amG_\classmap(\fk)$}{}.
We say that a Curtis-Tits amalgam has \nom{\em type $\Gamma_\classmap(\fk)$}{} if it is a homomorphic image of $\amG_\classmap(\fk)$. The Classification Theorem can thus be reformulated as follows.

\begin{namedtheorem}{Corollary}
Every Curtis-Tits amalgam with property (D) with connected simply-laced triangle free  {Dynkin} diagram over a field $\fk$ of size at least $4$ has type $\Gamma_\classmap(\fk)$ for some $\classmap$.
\end{namedtheorem}

\subsection{Curtis-Tits amalgams yielding Curtis-Tits groups}
For the remainder of the paper \nom{$\Gamma$}{}, \nom{$\classmap$}{}, and \nom{$\fk$}{} will be as in the Classification Theorem. 
We wish to decide which general Curtis-Tits amalgam of type $\Gamma_\classmap(\fk)$ have non-trivial universal completion.

\bde\label{dfn:relation 0}
Let \nom{$\sim_0$}{} be the equivalence relation on $I$ so that for distinct $i,j\in I$ we have $i\sim_0 j$ if and only if 
\begin{enumerate}
\item[(B1)] $i$ and $j$ are not adjacent, but have a common neighbor in $\Gamma$, and 
\item[(B2)] the neighbors of $i$ in $\Gamma$ coincide with the neighbors of $j$ in $\Gamma$.
\end{enumerate}
Now consider two equivalence relations $\sim$ and $\sim'$ on $I$. We say that \nom{$\sim$ {\dfn refines} $\sim'$}{} if and only if $i\sim j$ implies $i\sim' j$ for all $i,j\in I$; that is, each equivalence classes of $\sim$ is contained in an equivalence class of $\sim'$.
\ede 

\bde\label{dfn:amalgam relation} 
For any Curtis-Tits amalgam $\amA=\{\amgrpA_i,\amgrpA_{i,j},\ama_{i,j}\mid i,j\in I(\Gamma)\}$ of type $\Gamma_\classmap(\fk)$, let \nom{$\sim_\amA$}{} be the equivalence relation generated by the pairs $(i,j)\in I\times I$ such that for distinct $i,j\in $ we have $i\sim_\amA j$ if and only if 
\begin{enumerate}
\item[(B0)]  $\amgrpA_{i,j}\cong \SL_2(\fk)\times\SL_2(\fk)/\langle (z,z) \rangle$,
\end{enumerate}
where $\langle z\rangle =Z(\SL_2(\fk))$.
\ede

\bde\label{dfn:completion relation}
Next, suppose that the  Curtis-Tits amalgam $\amA=\{\amgrpA_i,\amgrpA_{i,j},\ama_{i,j}\mid i,j\in I(\Gamma)\}$ of type $\Gamma_\classmap(\fk)$ has a non-trivial completion $(\compA,\compa_\bullet)$.
Let  \nom{$\sim_\compA$}{} be the relation of the pairs $(i,j)\in I\times I$ such that for distinct $i,j\in $ we have $i\sim_\compA j$ if and only if 
\begin{enumerate}
\item[(B0)]  $\compa_{i,j}(\amgrpA_{i,j})\cong \SL_2(\fk)\times\SL_2(\fk)/\langle (z,z) \rangle$,
\end{enumerate}
where $\langle z\rangle =Z(\SL_2(\fk))$.
\ede
We shall prove in Proposition~\ref{prop:sim_R refines sim_0} that $\sim_\compA$ is in fact an equivalence relation, and that $\sim_\amA$ refines $\sim_\compA$.

Our first result allows one to decide which Curtis-Tits amalgams of type $\Gamma_\classmap(\fk)$ yield non-trivial Curtis-Tits groups. The proof is obtained in Subsection~\ref{subsec:non-trivial completion iff sim_A in B}.
\bmth\label{mth:non-trivial completion iff sim_A in B}
Let $\amA$ be a Curtis-Tits amalgam with property (D) of type $\Gamma_\classmap(\fk)$.
Then, $\amA$ has a non-trivial universal completion if and only if $\sim_\amA$ refines $\sim_0$.
\emth
Theorem~\ref{mth:non-trivial completion iff sim_A in B} says that $\sim_\amA$ must refine $\sim_0$.
For the standard Curtis-Tits amalgam $\amG_\classmap(\fk)$ of type $\Gamma_\classmap(\fk)$ the relation $\sim_{\amG_\classmap(\fk)}$ is the equality relation. Hence we obtain
\bmco\label{mth:non-trivial completions}\label{mco:non-trivial completions}
Any standard Curtis-Tits amalgam arising from the Classification Theorem has a non-trivial universal completion.
\emco
\bde\label{dfn:universal completion of standard amalgam}
{We shall denote \nom{$(\ucompG_\classmap(\fk),{\ucompg_\classmap(\fk)}_\bullet)$}{} as the universal completion of the standard Curtis-Tits amalgam $\amG_\classmap(\fk)$ of type $\Gamma_\classmap(\fk)$ and call $\ucompG_\classmap(\fk)$ the \nom{{\dfn simply-connected  Curtis-Tits group of type $\Gamma_\classmap(\fk)$}}{}. Since $\classmap$ and $\fk$ are fixed we often omit them from the notation.} 
The quotient of $\ucompG_\classmap(\fk)$ over its center will be called the 
\nom{\dfn adjoint Curtis-Tits group of type $\Gamma_\classmap(\fk)$}.\index{adjoint Curtis-Tits group of type $\Gamma_\classmap(\fk)$} From the proof of Theorem~\ref{mth:non-trivial completion iff sim_A in B} it will emerge that this group corresponds to the $\sim_0$ relation.
\ede

\subsection{Curtis-Tits amalgams injecting into their Curtis-Tits groups}
Curtis-Tits amalgams with a non-trivial completion do not necessarily inject into their universal completion. We now determine when they do.

Write \nom{$(\ucompG,\ucompg_\bullet)=(\ucompG_\classmap(\fk),{\ucompg_\classmap(\fk)}_\bullet)$}{}.
\bde\label{dfn:zij Z0 ZbfA}\label{dfn:closure}
For any amalgam $\amA$ of type $\Gamma_\classmap(\fk)$ with completion $(\compA,\compa_\bullet)$ consider the map \nom{$d^\compA\colon (\fk^*)^I\to \compA$}{} of Definition~\ref{dfn:di si}.
Define \index{$\kIz_{ij}$}\index{$\kIZ^0$}\index{$\kIZ^\compA$}
\begin{align*}
\kIz_{ij}{}&=(a_k)_{k\in I}\in (\fk^*)^I \mbox{ where } a_k=\begin{cases} -1 & \mbox{ if } k\in \{i,j\},\\ 1 &\mbox{ else}.\end{cases}\\
\kIZ^0&=\langle \kIz_{ij}\colon i\sim_0 j \mbox{ with }i,j\in I\mbox{ distinct}\rangle\le (\fk^*)^I.\\
\kIZ^\compA&=\ker d^\compA\cap \kIZ^0\le (\fk^*)^I.
\end{align*}
Let \nom{$\amG_1$}{} be the image of $\amG_\classmap(\fk)$ in $\ucompG$, let \nom{$\sim_1=\sim_{\amG_1}$}{} and let \nom{$\kIZ^1=\kIZ^{\ucompG}$}{}.
Using these groups, we define a map of equivalence relations on $I$ sending $\sim$ to the equivalence relation \nom{$\bar{\sim}$}{} generated by the pairs $(i,j)$ of distinct elements of $I$ such that  
\begin{align*}
\kIz_{ij}\in  \langle \kIZ^1, \kIz_{kl}\colon k\sim l\mbox{ with } k,l\in I \mbox{ distinct}\rangle.
\end{align*}
Clearly $\bar{\bar{\sim_\amA}}=\bar{\sim_\amA}$ so $\bar{\ \cdot\ }$ is a closure operator with respect to the refinement order.
\ede
We now have the following result. The proof is obtained in Subsection~\ref{subsec:non-trivial completion iff sim_A in B}.

\bmth\label{mth:amA injects iff sim_1 refines sim_A} 
Let $\amA$ be a Curtis-Tits amalgam with property (D) of type $\Gamma_\classmap(\fk)$ such that $\sim_\amA$ refines $\sim_0$.
Let $\bar{\amA}$ be the image of $\amA$ in its universal completion.
Then, $\sim_{\bar{\amA}}=\bar{\sim_\amA}$ hence $\amA$ injects into its universal completion if and only if $\sim_\amA=\bar{\sim_\amA}$.
\emth
To summarize Theorems~\ref{mth:non-trivial completion iff sim_A in B}~and~\ref{mth:amA injects iff sim_1 refines sim_A}, using properties (B1) and (B2) of Definition~\ref{dfn:relation 0} one can compute $\sim_0$ directly from $\Gamma$.
From this, one can determine which Curtis-Tits amalgams of type $\Gamma_\classmap(\fk)$ have a non-trivial universal completion. 
To find $\sim_1$ and the closure operator of Definition~\ref{dfn:closure}, one needs group theoretic properties of $\ucompG_\classmap(\fk)$. 
Once these are known we can find  the Curtis-Tits amalgams of type $\Gamma_\classmap(\fk)$ which inject into their universal completion as those corresponding to  the elements closed with respect to $\bar{\ \cdot \ }$ in the interval between $\sim_0$ and $\sim_1$ ordered by refinement.

\subsection{Construction of Curtis-Tits groups}
Let \nom{$\comop$}{} be the endomorphism of $(\fk^*)^I$ described in Definition~\ref{dfn:comop}.

In the non-orientable case, there exists a canonical orientable Curtis-Tits amalgam  \nom{$\amL$}{} over a $2$-sheeted covering \nom{$\La=(\hI,\hE)$}{} of $\Gamma$. Let {\nom{$\compL$}{} be the quotient of the universal completion of $\amL$ over its center. 
Then, the monodromy group of the covering of graphs $\La\to \Gamma$ induces an involution \nom{$\theta$}{} of $\compL$.  Let \nom{$L^\theta$}{} be the centralizer of $\theta$ in $\compL$ and let \nom{$\compG$}{} be its commutator subgroup. Let \nom{$\ucompD=\langle \ucompg_i(\amgrpD_i)\colon i\in I\rangle_\ucompG$}{}, with \nom{$\amgrpD_i$}{} as in Definition~\ref{dfn:property D}.

\bmth\label{thm:concrete completions}
First assume that $\Gamma_\classmap(\fk)$ is orientable, (that is $\im\classmap\le \Aut(\fk)$) and $|\fk|\ge 4$.
\begin{enumerate}
\item\label{orientable} The group $\ucompG$ is a central extension of  a group of Kac-Moody type over $\fk$ with diagram $\Gamma$. Moreover, $Z(\ucompG)=d(\ker \comop)$.
\end{enumerate}
Next, assume that $\Gamma_\classmap(\fk)$ is not orientable, (that is, $\im\classmap\not\le\Aut(\fk)$)  and that $|\fk|\ge 7$. With the notation above we have
\begin{enumerate}
\setcounter{enumi}{1}
\item  $(\compG,\compg_\bullet)$ is a non-trivial completion of $\amG$, for some $\compg_\bullet$.

\setcounter{enumi}{2}
\item  $\compLth/\compG$ is an elementary abelian $2$-group which can be explicitly computed using the operator $\comop$.

\item Suppose $\compLth=\compG$ and \nom{$\umap\colon \ucompG\to \compLth$}{} is the canonical map.
Then, $\ker\umap = Z(\ucompG)\cap\ucompD=d(\ker \comop)$.
 \end{enumerate}
\emth

Note that the first half of part 1 of Theorem~\ref{thm:concrete completions}, which was already accomplished in~\cite{BloHof2014b},  is similar to, but stronger than the Example following Theorem A in~\cite{Cap2007}. Indeed, we do not presuppose that the amalgam is isomorphic to the Curtis-Tits amalgam of a group of Kac-Moody type.
The proof of this part is reviewed in Subsection~\ref{subsec:orientable}
The second half, on the other hand is new.
The proofs of Part 2, Part 3, and Part 4 are obtained in Subsections~\ref{subsec:non-orientable CT groups},~\ref{subsec:parabolic amalgam}, and~\ref{sec:L=Gtheta} respectively.

We also record the following geometric result, arising from the proof of Theorem~\ref{thm:concrete completions}, which is of independent interest from a geometric group theory point of view. 
The proof of this theorem is obtained in Subsection~\ref{subsec:parabolic amalgam}.

\bmth\label{thm:G flag-transitive on Detheta}
Under the assumptions of part 2. of Theorem~\ref{thm:concrete completions}, $\ucompG$ acts flag-transitively on a connected and simply $2$-connected geometry $\De^\theta$.
\emth
The simple $2$-connectedness of the geometry $\De^\theta$ follows in a similar manner as in~\cite{BloHof2014a} (see Section~\ref{section:De theta}). The flag-transitivity of $\ucompG$ follows from Corollary~\ref{cor:BJ flag transitive}.

Finally, in Section~\ref{section:examples} we apply the methods developed in this paper to some concrete examples of orientable and non-orientable Curtis-Tits groups.

\section*{Acknowledgements}
Part of this paper was written as part of the project KaMCAM funded by the European Research Agency through a Horizon 2020 Marie-Sk\l odowska Curie fellowship (proposal number 661035).

The authors would also like to thank the excellent comments of the anonymous referee of an earlier version of the manuscript. They encouraged us to improve the presentation and strengthen some of the results.

\section{Definitions and properties of Curtis-Tits amalgams}\label{section:properties of CT amalgams}
\subsection{Amalgams and diagrams}\label{subsec:amalgams}
In this section we briefly recall the notion of a Curtis-Tits amalgam over a commutative field $\fk$ from~\cite{BloHof2014b} and prove that up to finite central extension one can restrict to universal Curtis-Tits amalgams. 
\medskip

\bde\label{dfn:amalgam}
An \nom{{\em amalgam} }{} over a poset $(\cP,\prec)$ is a collection $\amA=\{\amgrpA_x\mid x\in \cP\}$ of groups, together with a collection $\ama_\bullet=\{\ama_x^y\mid x\prec y, x,y\in \cP\}$ of monomorphisms $\ama_x^y\colon \amgrpA_x\into \amgrpA_y$, called {\em inclusion maps} such that whenever $x\prec y\prec z$, we have
 $\ama_x^z=\ama_y^z\after\ama_x^y$.
A {\em completion} of $\amA$ is a group $\compA$ together with a collection  $\compa_\bullet=\{\compa_x\mid x\in \cP\}$ of homomorphisms $\compa_x\colon \amgrpA_x\to \compA$, whose images generate $\compA$, such that for any $x, y\in \cP$ with $x\prec y$ we have 
$\compa_y\after\compa_x^y=\compa_x$.
The amalgam $\amA$ is {\em non-collapsing} if it has a non-trivial completion.
A completion $(\ucompA,\ucompa_\bullet)$ is called {\em universal} if for any completion $(\compA,\compa_\bullet)$ there is a unique surjective group homomorphism $\umap\colon \ucompA\to \compA$ such that $\compa_\bullet=\umap\after\ucompa_\bullet$. 
\ede
\bre
Throughout the paper we will adhere to the font conventions in Definition~\ref{dfn:amalgam}, that is, amalgams are in calligraphic font, their groups and connecting maps are in boldface roman, 
 their completion maps are in greek and their image is denoted in regular math font.
 Finally, universal completions use the notation of ordinary completions with a tilde on top.
 \ere
We record the following lemma, which will be applied in the proof of Lemma~\ref{lem:tK=ker phi} (for another example of its use see the proof of Proposition 3.12 of~\cite{Cap2007}).
\ble\label{lem:locally generated kernel}
Let $\eta\colon A\to B$ be a surjective group homomorphism.
Suppose $\amA=\{\compA_x \mid x \in \cP\}$ is an amalgam of subgroups of $\compA$, with connecting maps given by inclusion,  generating $\compA$.
Let $\amB=\{\compB_x=\eta(\compA_x)\mid x\in \cP\}$ be the amalgam corresponding to $\amA$ via $\eta$ and suppose that 
 $\compB$ is the universal completion of $\amB$.
For each $x\in \cP$, let $K_x=\ker \eta|_{\compA_x}$, 
then,  $\ker\eta$ is the normal closure of the subgroup generated by all $K_x$.
\ele
\bpf
Let $K$ be the normal closure of the subgroup generated by all $K_x$ ($x\in \cP$).
Then, by definition $K\le \ker \eta$.
Moreover, for any $x\in \cP$ we have
 $K_x\le K\cap \compA_x\le \ker \eta\cap \compA_x=K_x$.
Therefore, the canonical map $\eta'\colon \compA\to \compA/K$ restricts to $\compA_x\to \compA_x/K_x\stackrel{\eta_x}{\cong} \compB_x$.
By universality of $\compB$ there exists a map 
 $\compB\to \compA/K$ extending this isomorphism.
On the other hand, since $K\le \ker \eta$, there is also a homomorphism $\compA/K\to \compB$ induced by $\eta$. 
The composition $\compB\to \compA/K\to \compB$ fixes the elements of $\amB$ elementwise so it must be the identity mapping.
\epf

\bde\label{dfn:diagram}
For the purposes of this paper, a \nom{simply-laced diagram}{} is an undirected graph {$(I,E)$}{} with finite vertex set $I\sbe \NN$ and edge set $E$ without circuits of length $1$ or $2$. All simply-laced diagrams in this paper are {connected and} triangle-free, that is, they also have no circuits of length $3$.
\ede
Let \nom{$\Gamma=(I,E)$}{}  denote a {connected} triangle-free simply-laced diagram with $3\le n=|I|<\infty$ nodes.\idx{$n$}
Also fix a field \nom{$\fk$}{} of order at least $4$.

\paragraph{Indexing convention}
Throughout the paper we shall adopt the following indexing conventions. Indices from $I$ shall be taken modulo $n=|I|$. For any $i\in I$, we set $(i)=I-\{i\}$.
Also subsets of $I$ of cardinality $1$ or $2$ appearing in subscripts are written without set-brackets.

\bde\label{dfn:CT structure}
Let $\cP=\{J\mid \emptyset\ne J\sbe I \mbox{ with }|J|\le 2\}$ and 
 $\prec$ denoting inclusion.
A  {\em Curtis-Tits amalgam with diagram  $\Gamma$ over $\fk$}  is an amalgam 
\nom{$\amG=\{\amgrpG_{i},\amgrpG_{i, j}, \amg_{i,j} \mid  i, j \in I\}$}{} 
over $\cP$, where, for every $i,j\in I$, we write $\amg_{i,j}=\amg_{\{i\}}^{\{i,j\}}$.
Note that, due to our subscript conventions, we write  $\amgrpG_{i}=\amgrpG_{\{i\}}$  and $\amgrpG_{i,j}=\amgrpG_{\{i,j\}}$, 

where 
\begin{enumerate}
  \item[(SCT1)]  for any vertex $i$, we set $ \amgrpG_i = \SL_2(\fk)$ and for each pair $i,j \in I$,
  $$\amgrpG_{i,j}\cong\begin{cases}
  \SL_3(\fk) & \mbox{if} \ \{i,j\}\in E \\ \amgrpG_i\circ\amgrpG_j & \mbox{else}\end{cases};$$
  \item[(SCT2)]
    For $\{i,j\}\in E$, $\amg_{i,j}(\amgrpG_i)$ and $\amg_{j,i}(\amgrpG_j)$ form a standard pair for $\amgrpG_{i,j}$ in the sense of~\cite{BloHof2014b}, whereas for all other pairs $(i,j)$, $\amg_{i,j}$ is the natural inclusion of $\amgrpG_i$ in $\amgrpG_i\circ\amgrpG_j$.
\end{enumerate}

{
Here, \nom{$\circ$}{ denotes central product}. Note that there are two cases:
\begin{align*}
\amgrpG_i\circ \amgrpG_j \cong\begin{cases}
\amgrpG_i\times \amgrpG_j,\\
\amgrpG_i\times \amgrpG_j/\langle z,z\rangle, 
\end{cases}
\end{align*}}
{where $\langle z\rangle=Z(\SL_2(\fk))$.}
{We call the amalgam $\amG$ {\em universal} if the first case occurs for any $\{i,j\}\not\in E$. This happens for instance if $\Char(\fk)=2$.}\idx{universal amalgam}
\ede

We recall Property (D) from~\cite{BloHof2014b}.
\bde\label{dfn:property D}\nom{{(\rm property (D))}}{}
For any $i,j\in I$ with $\{i,j\}\in E$, let 
\begin{align*}
\amgrpD_i^j=N_{\amgrpG_{i,j}}(\amg_{j,i}(\amgrpG_j))\cap \amg_{i,j}(\amgrpG_i)
\end{align*}
Once checks that \nom{$\amgrpD_i^j$}{} is the only torus in $\amg_{i,j}(\amgrpG_i)$ normalized by $\amgrpD_j^i$ (see~\cite{BloHof2014b}).
We now say that $\amG$ {\dfn has property (D)} if for any pair of edges $\{i,j\}$ and $\{i,k\}$ in $E$ the map $\amg_{i,k}\after\amg_{i,j}^{-1}\colon \amgrpD_i^j\to \amgrpD_i^k$ is an isomorphism.
{That is, we can select tori \nom{$\amgrpD_i\le \amgrpG_i$}{} such that for all edges $\{i,j\}$ we have
     $\amg_{i,j}(\amgrpD_i)=\amgrpD_i^j$.}
\ede
\subsection{Spherical Curtis-Tits amalgams}\label{subsec:spherical Curtis-Tits}
In the spherical case we have the following corollary to a version of the Curtis-Tits theorem due to Timmesfeld~\cite[Theorem 1]{Tim1998}:
\bth\label{thm:Tim1998}
Let $\amG$ be a Curtis-Tits amalgam (universal or not)  over a commutative field $\fk$, having spherical diagram $\Gamma$. Then any non-trivial  completion $\compG$ of $\amG$ is a perfect central extension of the adjoint Chevalley group of type $\Gamma$ and a central quotient of the universal Chevalley group of type $\Gamma$.
\eth
\bpf
The condition (H1) preceding Theorem 1 of~\cite{Tim1998} is satisfied.
For each $i,j\in I$, let $X_i=\amgrpG_i$, and $H_i=\amgrpD_i$ (the standard torus $\amgrpD_i$ arising from Property (D)).
It follows that, for each $i,j\in I$ we have $X_{i,j}=\amgrpG_{i,j}$.

We verify the conditions in (H1):
(1) We have $X_i\cong\SL_2(\fk)$, which is a central extension of $\PSL_2(\fk)$.

(3) Since the images of $X_i$ and $X_j$ form a standard pair in $X_{i,j}\cong \SL_3(\fk)$, whenever $\{i,j\}$ is an edge of $\Gamma$, and $X_{i,j}$ is the central product of the images of $X_i$ and $X_j$ whenever $\{i,j\}$ is not an edge in $\Gamma$, this condition is satisfied.

(2) Property (D) ensures that this condition is satisfied as well.
We also check the rest of the conditions of Theorem 1 of loc. cit.. Namely, if $|\fk|=4$ {and $\{i,j\}$ is an edge of $\Gamma$, then $\amgrpG_{i,j}=\SL_3(4)$, which has 
  $|Z(\SL_3(\fk))|=3$.}
The result follows.
\qed 
\bre\label{rem:universal and adjoint groups}
From~\cite{GorLyoSol1998} we have the following identifications.
\begin{enumerate}
\item The universal and adjoint Chevalley groups of type $A_n$ are $\SL_{n+1}(\fk)$ and $\PSL_{n+1}(\fk)$.

\item
The universal and adjoint Chevalley groups of type $D_n$ over $\fk$ are the spin group $\Spin_{2n}^+(\fk)$ and $\POmega^+(\fk)$, the derived subgroup of the group of linear isometries of a non-degenerate quadratic form of Witt index $n$, modulo the center of this derived subgroup. 
\item We have $\Spin_{2n}^+(\fk)\cong 2^2 \POmega^+_{2n}(\fk)$ if $n$ is even and $\Spin_{2n}^+(\fk)\cong 4 \POmega^+_{2n}(\fk)$ if $n$ is odd~\cite{Wil2009}.
 
\end{enumerate}
\ere

\subsection{General and Universal Curtis-Tits amalgams}\label{subsec:general universal}
Note that every Curtis-Tits amalgam is a quotient of a unique universal one that we call $\amG$.
This means in particular that any completion of a Curtis-Tits amalgam is a completion of the corresponding universal amalgam.
We shall now investigate the relation between the universal completions of universal and other Curtis-Tits amalgams.

Let $\amG$ be a universal Curtis-Tits amalgam with {connected simply-laced} triangle-free diagram  $\Gamma$ over $\fk$. We shall only be interested in amalgams that admit a non-trivial completion $(\compA,\compa_\bullet)$. From~\cite{BloHof2014b} it then follows that $\amG$ has property (D).
Using connectedness of $\Gamma$ and the fact that conjugates of the vertex groups generate the edge groups, one shows that in fact all maps $\compa_J$ ($J\sbe I$ with $0<|J|\le 2$) are non-trivial.

\paragraph{Notation}
Let \nom{$(\compA,\compa_\bullet)$}{} be a completion of $\amG$ such that for each $i,j\in I$, the maps $\compa_i$ and $\compa_{i,j}$ are non-trivial. We discuss the structure of the image of $\amG$ in $\compA$. We show that the possibilities are only limited by the structure of $\Gamma$.

Given a subset \nom{$S\sbe I$}{}, let  \nom{$\Gamma_S$}{} be the full subgraph of $\Gamma$ supported by the node set $S$.
We define the subamalgam \nom{$\amG_S=\{\amgrpG_J\mid J\sbe S\mbox{ with }0<|J|\le 2\}$}{} of $\amG$, with connecting maps induced by $\amg_\bullet$. Let \nom{$\compA_S$}{} be the subgroup of $\compA$ generated by the images of the elements of $\amG_S$ under $\compa_\bullet$, let \nom{$(\ucompG_S,\ucompg_{S;\bullet})$}{} be the universal completion of $\amG_S$, and let \nom{$\umap_S\colon \ucompG_S\to \compA_S$}{} be the universal map.
For $S=I$, we set $\ucompG=\ucompG_S$ {and then $\umap=\umap_S$}.

\paragraph{The maps $\compa_i$ and $\compa_{i,j}$}
\ble\label{lem:SL2 SL3}\label{lem:vertex and edge group images}
The map $\compa_i$ is injective for all $i\in I$ and $\compa_{i,j}$ is injective for all $\{i,j\}\in E$.
\ele
\bpf 
Since {the diagram $\Gamma$} is connected and has at least three nodes, but no triangles, every $i\in I$ lies on an edge and any edge is part of a subdiagram of type   $A_3$. 
Without loss of generality let us assume that $i,j,k\in I$
 are such that $j$ is adjacent to both $i$ and $k$ and let $S=\{i,j,k\}$.
{By  Theorem~\ref{thm:Tim1998}~and~Remark~\ref{rem:universal and adjoint groups} }we see that  the universal completion $\ucompG_S$ of the amalgam $\amG_S$ is a central extension of the group $\compA_S$  isomorphic to $\SL_4(\fk)$.   
For each $s\in S$, let $\ucompU_{s}^\pm$ be the root groups mapping isomorphically to the standard root groups \nom{$\compU_s^{\pm}$}{} of $\compA_s$ under $\umap_S$.
Then, $\SL_2(\fk)\cong \langle \ucompU_s^+,\ucompU_s^-\rangle$ does not intersect the center of $\ucompG_s$ and hence maps isomorphically to $\compA_s$. For the same reason $\SL_3(\fk)\cong \langle \ucompU_i^\pm , \ucompU_j^\pm\rangle$ maps isomorphically to $\compA_{i,j}$ and the same holds for $j$ and $k$. 
\epf

\medskip
Next we consider the image of $\compa_{i,j}$, when $\{i,j\}\not\in E(\Gamma)$.
\bde\label{dfn:amz}
If $\Char \fk$ is odd, then by \nom{$\amz_i$}{} we denote the unique central involution of $\amgrpG_i\cong \SL_2(\fk)$.
By Lemma~\ref{lem:vertex and edge group images}, the  image of $\amz_i$ is non-trivial in any completion of $\amG_S$, where $\Gamma_S$ is connected of size at least $2$.
We denote \nom{$\ucompz_i=\ucompg_{S;i}(\amz_i)\in \ucompG_S$}{}.
\ede

\ble\label{lem:general non-edge kernel}
Let $\{i,j\}\not\in E$.
If $\Char \fk$ is even then $\compa_{i,j}$ is injective.
Otherwise $\ker \compa_{i,j}\le \langle (\amz_i,\amz_j)\rangle\le \amgrpG_i\times \amgrpG_j$.
\ele
\bpf
By Lemma~\ref{lem:vertex and edge group images}, $\compa_i$ and $\compa_j$ are injective and their images commute.
{As $\compa_{i,j}$ extends both maps, the kernel must be central in $\amgrpG_i\times \amgrpG_j$.}
\epf
\paragraph{The equivalence relation $\sim$ of bad pairs}
\bde
With the notation of Lemma~\ref{lem:general non-edge kernel} we call $\{i,j\}$ a \nom{ {\em good pair}}{} if $\ker\compa_{i,j}$ is trivial and a \nom{ {\em bad pair}}{} otherwise.
In case $\{i,j\}$ is a bad pair we write \nom{$i\sim_\compA j$}{}.
\ede

The following is verified easily inside $\SL_4(\fk)$.
\ble\label{lem:A3}
If $S=\{i,k,j\}$ is such that $\Gamma_S$ has diagram $A_3$, and $k$ is a neighbor of $i$ and $j$, then 
 $\compa_{i,j}(\amz_i\amz_j)\in Z(\compA_S)$.
\ele
 
\ble\label{lem:A1 x A2}
If $\Gamma_S$ has diagram $A_1\times A_2$, then, $(S,\sim_\compA)$ has no edges.

\ele
\bpf
Let $S=\{i,j,k\}$ with $i$ having no neighbor in the edge $\{j,k\}$.
If $\{i,j\}$ is a bad pair, then  {$\compa_{i,j}(\amz_i,\amz_j)=e\in G$}.
This means that the non-trivial elements $\compa_i(\amz_i)$ and $\compa_j(\amz_j)$ are equal.
However, in $\amG$, we see that $\amz_i$ commutes with all of $\amgrpG_k$, but $\amz_j$ does not. Note that by Lemma~\ref{lem:vertex and edge group images} the $\compa_{j,k}$ is injective so $\compa_j(\amz_j)$ does not commute with all of $\compa_k(\amgrpG_k)$, whereas $\compa_i(\amz_i)$ does, a contradiction.
\epf

\ble\label{lem:D4 lemma}
Let $S=\{i,j,k,l\}$ be such that $\Gamma_S$ has diagram $D_4$ with central node $k$.
Then, $(\{i,j,l\},\sim_\compA)$ has one or three edges.
\ele
\bpf
For each $s\in \{i,j,k,l\}$ let $\ucompz_s=\ucompg_{S;s}(\amz_s)\in \ucompG_S$.
This is non-trivial by Lemma~\ref{lem:vertex and edge group images}.

{From Theorem~\ref{thm:Tim1998}~and~Remark~\ref{rem:universal and adjoint groups} we know that $\ucompG_S$ is the universal Chevalley group $\Spin^+_8(\fk)$  having $\compA_S$ as a central quotient.
Moreover, from that remark we also have that $C_2\times C_2\cong \{1,\ucompz_i\ucompz_j,\ucompz_j\ucompz_l,\ucompz_i\ucompz_l\}= Z(\ucompG_S)$. }The conclusion follows.
\epf

\ble\label{lem:bad pair partition}
The graph $(I,\sim_\compA)$ is a disjoint union of complete graphs.
 \ele
\bpf
Suppose $i\sim_\compA j$.
By Lemma~\ref{lem:A1 x A2}, and connectedness of $\Gamma$, $i$ and $j$ have a common neighbor $k$ in $\Gamma$. 
Now if $i\sim_\compA l$, then by Lemma~\ref{lem:A1 x A2}, $l$ is a neighbor of $k$ also.
By Lemma~\ref{lem:D4 lemma} since $j\sim_\compA i\sim_\compA l$, also $j\sim_\compA l$.
\epf

\medskip
Let \nom{$\cB=\{B_1, \ldots, B_m\}$}{} be the set of connected components of $(I,\sim_\compA)$.  Lemma~\ref{lem:A1 x A2} implies the following.

\bco\label{cor:sim_cB}
If  $i\in B_s$ and $j\in B_t$ are adjacent in $\Gamma$, then $\Gamma_{B_s\cup B_t}$ is a complete bipartite graph.
\eco
\bpr\label{prop:sim_R refines sim_0}
Suppose that the Curtis-Tits amalgam $\amA$ of type $\Gamma_\classmap(\fk)$ has a non-trivial completion $(\compA,\compa_\bullet)$. Then, $\sim_\compA$ is an equivalence relation. Moreover, $\sim_\amA$ refines $\sim_\compA$ and $\sim_\compA$ refines $\sim_0$.
\epr
\bpf
The fact that $\sim_\compA$ is an equivalence relation is Lemma~\ref{lem:bad pair partition}.
Since the map $\compa_\bullet$ can only introduce more bad pairs, clearly $\sim_\amA$ refines $\sim_\compA$.
It follows from Corollary~\ref{cor:sim_cB} that $\sim_\compA$ refines $\sim_0$.
\epf

\paragraph{Isogeny of Curtis-Tits amalgams}
Define a graph 
 \nom{$(\cB,\sim_{\cB})$}{}, where $B_s\sim_{\cB} B_t$ if $\Gamma_{B_s\cup B_t}$ is a complete bipartite graph.
\ble\label{lem:cB graph}
The graph $(\cB,\sim_{\cB})$ is connected.
\ele
\bpf
Connectedness follows from Corollary~\ref{cor:sim_cB}  and connectedness of $\Gamma$.
\epf

\medskip

\bco\label{cor:bad graph characterization}
Suppose $\cB=\{B_i\mid i=1,\ldots,m\}$.
Then, $I=\biguplus_{i=1}^m B_i$ and for any selection $x_i\in B_i$ ($i=1,\ldots,m$)
 we have that $\{x_i,x_j\}\in E$ if and only if $B_i\sim_\cB B_j$.
\eco

\bco\label{cor:center}
There exists a central elementary abelian $2$-subgroup $\ucompZ$ of $\ucompG$  of size at most $2^{|I|-|\cB|}$ such that $\ucompG/\ucompZ$ contains a copy of the image of $\amG$ in $\compA$.  In fact $\ucompZ$ is a central subgroup contained in the diagonal group $\ucompD$.
\eco
\bpf
For each $B\in \cB$, let 
 $\ucompZ_B$ be the subgroup generated by all products $\ucompz_i\ucompz_j\in \ucompG=\ucompG_I$, where $i,j\in B$ (so in case $|B|=1$ $\ucompZ_B$ is trivial). Thus $\ucompZ_B$ is an elementary abelian group of size at most $2^{|B|-1}$.
 By Lemma~\ref{lem:A3} and since $\amz_s\in \amgrpD_s$, for all $s\in I$, we have $\ucompZ_B\le Z(\ucompG)\cap \ucompD$.  
Take $\ucompZ=\langle \ucompZ_B\mid B\in \cB\rangle$. Then $\ucompZ\le Z(\ucompG)\cap \ucompD$ is an elementary abelian $2$-group of size at most 
 $2^{|I|-|\cB|}$. 
Consider the map $\umap\colon \ucompG\to \compA$.
From Lemma~\ref{lem:vertex and edge group images} we know that $\ker\umap$ intersects 
 none of the groups $\ucompG_i$ and $\ucompG_{i,j}$ for $i\in I$ and $\{i,j\}\in E$.
Moreover, $i\sim_\compA j$ if and only if $\ucompz_i\ucompz_j\in \ker\umap$.
Thus, $\ucompZ\le \ker \umap$ so that $\umap$ factors as $\ucompG\to \ucompG/\ucompZ\stackrel{\underline{\umap}}\to \compA$.
Now consider distinct $i,j\in I$ such that $\{i,j\}\not\in E$. Then by Lemma~\ref{lem:A1 x A2}
 and since $\Gamma$ is triangle free and connected, $\alpha_{i,j}(\amz_i\amz_j)=1$ implies $i\sim_\compA j$ and conversely. It follows that the map $\underline{\umap}$ is the identity on the image of $\amG$.
\epf

\bre
Note that $\ker\underline{\umap}$ in the proof of Corollary~\ref{cor:center} does not have to be trivial.
\ere
\medskip

Given $(\compA,\compa_\bullet)$, we shall denote the central elementary abelian $2$-group of Corollary~\ref{cor:center} by \nom{$\ucompZ_\compA$}{}.

\bde\label{dfn:isogenic}
{We call two groups $A$ and $B$ \nom{ {\em isogenic} }{}if there is a group $C$ which is a finite central quotient of both $A$ and $B$.
}
\ede

\bmth\label{thm:isogenic completions}
Given a Curtis-Tits amalgam $\amA$ {of type $\Gamma_\classmap(\fk)$} with non-trivial completion, then its universal completion $(\ucompA,\ucompa_\bullet)$ is isogenic to the universal completion of the corresponding universal Curtis-Tits amalgam of {type $\Gamma_\classmap(\fk)$}.
More precisely, we have $\ucompA=\ucompG/\ucompZ_\ucompA$.
\emth
\bpf
Let $\amG$ be the universal Curtis-Tits amalgam  corresponding to $\amA$. This means that there is a map of amalgams $\ammapf_\bullet\colon\amG\to \amA$.
Let $(\ucompG,\ucompg_\bullet)$ be the universal completion of $\amG$ and let $(\ucompA,\ucompa_\bullet)$ be the universal completion of $\amA$. 
Then $\ammapf_\bullet$ induces a map $\ucompmapf\colon \ucompG\to \ucompA$.
Thus, $(\ucompA,\ucompmapf\after \ucompg_\bullet)$ is a completion of $\amG$.
We now apply Corollary~\ref{cor:center}. First note that the image of $\amG$ in $\ucompA$ is also the image of $\amA$ in $\ucompA$.
Since $\ucompG/\ucompZ_\ucompA$ contains a copy of $\amA$, $(\ucompG/\ucompZ_\ucompA,\underline{\ucompa}_\bullet)$ is a completion of $\amA$ for some map $\underline{\ucompa}_\bullet$. The proof also establishes a map $\underline{\ammapf}\colon \ucompG/\ucompZ_\ucompA\to \ucompA$ (called $\underline{\umap}$ there). By universality of $(\ucompA,\ucompa_\bullet)$ as a completion of $\amA$, we must have $\ucompA=\ucompG/\ucompZ_\ucompA$.
Since $\ucompZ_\ucompA$ is finite, we are done.
\epf

\bre
Note that in~\cite[Proposition 3.12]{Cap2007} it is proved that given a {Kac-Moody group} $G$ of $2$-spherical type $\Gamma$, and given the groups $G^\rsc$, which is the simply-connected Kac-Moody group of the same type, and $\check{G}$ the universal completion of the Curtis-Tits amalgam contained in $G$, then, the canonical map $\umap\colon G^\rsc\to \check{G}$ has finite central kernel.
In the proof of Corollary~\ref{cor:center}, $\ucompG$ and $\compA$ play the roles of $G^\rsc$ and $\check{G}$ respectively, so in fact the kernel of $\umap$ is a central $2$-group of size at most $2^{|I|-|\cB|}$.
\ere

\subsection{{Universal Curtis-Tits amalgams; standard Curtis-Tits amalgams}}\label{subsec:universal CT amalgams}
Because of Theorem~\ref{thm:isogenic completions} we shall now restrict our attention to universal Curtis-Tits amalgams and their completions.
These universal Curtis-Tits amalgams were classified in~\cite{BloHof2014b} (see Classification Theorem in the Introduction).

In this subsection we will make this correspondence precise and define one standard universal Curtis-Tits amalgam for each $\classmap\colon \Pig\to \Aut(\fk)\times\langle \tau\rangle$.
Before we do this, we specify an action of the group \nom{$\Aut(\fk)\times\langle \tau\rangle$}{} (with $\tau$ of order $2$) on $\SL_2(\fk)$.
We let $\alpha\in \Aut(\fk)$ act entry-wise on the matrix $A\in \SL_2(\fk)$ and let $\tau$ act by sending each $A\in \SL_2(\fk)$ to its transpose inverse ${}^tA^{-1}$ with respect to the standard basis. Note that $\tau$ acts as an inner automorphism. 

Fix a spanning tree \nom{$T=(I,E(T))$}{} for $\Gamma$ rooted at $0\in I$ and suppose that for any $e\in E(T)$,{ the vertex of $e$ nearest to $0$ (in $T$) has the smallest label}.
Given $e=\{i,j\}\in E-E(T)$ with $i<j$, let \nom{$\gamma(e)$}{} be the (the homotopy class of) the unique minimal  loop of the graph $(I,E(T)\cup \{e\})$ directed such that $e$ is traversed from $i$ to $j$.
This establishes an isomorphism of free groups $F(E-E(T))\to \Pig$.
Thus the homomorphisms $\classmap$ of the Classification Theorem correspond bijectively to {set maps $E-E(T)\to \Aut(\fk)\times\langle\tau\rangle$ given by $e\mapsto\classmap(\gamma(e))$; we shall denote this set map also by \nom{$\classmap$}{}.}

\bde\label{dfn:standard CT structure}
Let $\cP=\{J\mid \emptyset\ne J\sbe I \mbox{ with }|J|\le 2\}$ and 
 $\prec$ denoting inclusion.
Given a set-map $\classmap\colon E-E(T)\to \Aut(\fk)\times\langle\tau\rangle$, the \nom{{\em standard Curtis-Tits amalgam of type $\Gamma_\classmap(\fk)$}}{} is the amalgam 
\nom{$\amG^\classmap(\fk)=\{\amgrpG_{i},\amgrpG_{i, j}, \amg_{i,j} \mid  i, j \in I\}$}{} 
over $\cP$, where, for every $i,j\in I$, we write 
 \begin{align*}
\amgrpG_{i}=\amgrpG_{\{i\}}, && \amgrpG_{i,j}=\amgrpG_{\{i,j\}}, && \mbox{ and } \amg_{i,j}=\amg_{\{i\}}^{\{i,j\}}
\end{align*}
where 
\begin{enumerate}
  \item[(SCT1)]  for any vertex $i$, we set $ \amgrpG_i = \SL_2(\fk)$ and for each pair $i,j \in I$,
  $$\amgrpG_{i,j}\cong\begin{cases}
  \SL_3(\fk) & \mbox{if} \ \{i,j\}\in E \\ \amgrpG_i\times\amgrpG_j & \mbox{else}\end{cases};$$
  \item[(SCT2)]
  For $\{i,j\}\in E$ with $i<j$ we have 
  $$ 
  \begin{array}{cc}
  \begin{array}{rl}
\amg_{i,j}\colon \amgrpG_{i}&\to \amgrpG_{i,j}\\ 
 & \\
A & \mapsto  \left(\begin{array}{@{}cc@{}}A & 0 \\0 & 1\end{array}\right)
\end{array}
&
 \begin{array}{rl}
\amg_{j,i}\colon \amgrpG_{j}&\to \amgrpG_{j,i}\\
 & \\
A & \mapsto  \left(\begin{array}{@{}cc@{}}1 & 0 \\0 & A^{\classmap_{j,i}}\end{array}\right)
\end{array}
\end{array},$$

whereas for all other pairs $(i,j)$, $\amg_{i,j}$ is the natural inclusion of $\amgrpG_i$ in $\amgrpG_i\times\amgrpG_j$.

\end{enumerate}
Here 
$$
\classmap_{j,i}=
\begin{cases}
\classmap(\{j,i\}) & \mbox{ if }\{j,i\}\in E-E(T)\mbox{ and } i<j\\
1 & \mbox{ else. }
\end{cases} 
$$
\idx{$\classmap_{\{i,j\}}$}
We call $\amG^\classmap(\fk)$ {\em orientable} if $\im\classmap\le \Aut(\fk)$.
\ede
This notation directly generalizes that of~\cite{BloHof2014a}.
By~\cite{BloHof2014b} every Curtis-Tits amalgam over $\fk$ with {connected} diagram $\Gamma$ is of the form $\amG^\classmap(\fk)$ for some $\classmap$.

\subsection{Combinatorics of centers in Curtis-Tits groups}\label{subsec:centers in CT groups}
Let $(\compG, \compg_\bullet) $ be a completion of a universal Curtis-Tits amalgam  $\amG=\{\amgrpG_{i},\amgrpG_{i, j}, \amg_{i,j} \mid  i, j \in I\}$ {of type $\Gamma_\classmap(\fk)$ for some $\classmap$.}

\idx{$d_i^\compG$}\idx{$s_i$}\idx{$d^\compG$}
\bde\label{dfn:di si}
For any $i\in I$ and $a\in \fk^*$, let 
\begin{align*}
d_i^\compG(a)&=\compg_i\left(\begin{array}{@{}cc@{}} a & 0 \\ 0 & a^{-1} \end{array}\right),\\
s_i&=\compg_i\left(\begin{array}{@{}cc@{}}0 & 1\\ -1 & 0 \end{array}\right).
\end{align*}
Let \nom{$\compD$}{} be the subgroup of $\compG$ generated by the  tori \nom{$\compD_i=\im d_i$}{}.
Then \nom{$d^\compG\colon (\fk^*)^I\to \compD$}{} given by $d^\compG((a_i)_{i\in I})=\prod d_i^\compG(a_i)$ is a surjective homomorphism.

Note that by Lemma~\ref{lem:vertex and edge group images}, $\compg_i$ is injective for all $i\in I$ and, if {$\{i,j\}\in E$},  then also $\compg_{i,j}$ is injective.
Moreover, note that in this case
 $\compg_i=\compg_{i,j}\after \amg_{i,j}$.
 Therefore any computation in the subgroup 
  $\langle d_i^\compG(x), s_i, d_j(y), s_j\mid x, y\in \fk^*\rangle$ can be done entirely inside of $\amgrpG_{i,j}$.
If $\compG$ is clear from the context, we shall often omit the superscript $\compG$.
\ede

We will investigate the group $Z(\compG)\cap \compD$. 
To this end we have the following 
\bde\label{dfn:comop}
For any $i,j\in I$ with $i\sim j$, let 
 \nom{$\rho_{j,i}=\classmap_{i,j}^{-1}\classmap_{j,i}\in \Aut(\fk)\times\langle \tau\rangle\le \Aut(\fk^*)$}{}.
 Here $\tau$ embeds as the map $x\mapsto x^{-1}$.
We define the map $\comop=\comop_\amG=(\comop_i)_{i\in I}\colon (\fk^*)^I\to (\fk^*)^I$\idx{$\comop$},  by setting, for $a=(a_i)_{i\in I}\in (\fk^*)^I$
\begin{align*}
\comop_i(a)=a_i^{-2}\prod_{j\sim i} \rho_{j,i}(a_j).
\end{align*}
\ede 

We can now describe the action of $d(a)$ on $\amgrpG_i$.
\ble\label{lem:torus action}
Suppose that $d(a)=d^\compG(a)=\prod_{j\in I} d_j(a_j)$ and that 
 $$g=\compg_i\begin{pmatrix} x & y \\ w & z \end{pmatrix}\in \compg_i(\amgrpG_i).$$
Then, 
 $$(d(a))^{-1} g d(a) = \compg_i\begin{pmatrix} x & y k \\ k^{-1}w & z \end{pmatrix},$$
 where $k=\comop_i((a_j)_{j\in I})$.
\ele
\bpf
First note $d_j(a_j)=\compg_j\begin{pmatrix} a_j & 0 \\ 0 & a_j^{-1} \end{pmatrix}=\compg_{j,i}\left(\amg_{j,i}\begin{pmatrix} a_j & 0 \\ 0 & a_j^{-1} \end{pmatrix}\right)$.
This equals
\begin{align*}
\compg_{j,i}\left(\begin{pmatrix} \classmap_{j,i}(a_j) & 0 & 0  \\ 0 & \classmap_{j,i}(a_j)^{-1}  & 0 \\ 0 & 0 & 1 \end{pmatrix}\right) & \mbox{ if } j<i,\\
\compg_{j,i}\left(\begin{pmatrix} 1  & 0 & 0\\ 
0 & \classmap_{j,i}(a_j) & 0  \\ 0 & 0 & \classmap_{j,i}(a_j)^{-1} \end{pmatrix}\right) & \mbox{ if } j>i.
\end{align*}
In both cases conjugation gives
$$(d_j(a_j))^{-1} g d_j(a_j) = \compg_i\begin{pmatrix} x & y \rho_{j,i}(a_j) \\ \rho_{j,i}(a_j)^{-1} w & z \end{pmatrix}.$$
Finally, note that 
 $$(d_i(a_i))^{-1} g d_i(a_i) = \compg_i\begin{pmatrix} x & y a_i^{-2} \\ a_i^2w & z \end{pmatrix}.$$
The result follows since the actions of all $d_j(a_j)$ ($j\in I$) commute.
\epf
 
 \medskip

The following is an immediate consequence of Lemma~\ref{lem:torus action}.
\bth\label{prop:split null L}\label{thm:comop characterizes center}
Let $a=(a_i)_{i\in I}\in (\fk^*)^I$.
Then, $d^\compG(a)\in Z(\compG)$ if and only if $a\in \ker\comop$.
\eth

\bre
Note that $\comop$ does not depend on $\compG$, but $d^\compG$ does.
\ere

\bth\label{thm:KM center}
Suppose that $\compG$ is a {group of Kac-Moody type} and $\amG$ is the Curtis-Tits amalgam arising from the action on the twin-building \nom{$\De$}{} associated to its twin BN-pair.
Then the following are equivalent for $a\in (\fk^*)^I$.
\begin{enumerate}
\item $d^\compG(a)\in Z(\compG)$, 
\item $\comop(a)=(1)^I$,
\item $d^\compG(a)$ acts trivially on $\De$.
\end{enumerate}
\eth
\bpf
By Theorem~\ref{thm:comop characterizes center}, 1.~and 2.~are equivalent and by Lemma 3.4 of~\cite{Cap2007} 3.~and 1.~are equivalent.
\epf

\bre
Note that in the case when $G$ is a group of Lie type, this result is equivalent to the one of Steinberg (see for example section 12.1 of~\cite{Ca1972}) but is more elementary as it does not make use of weight lattices. 
Moreover, this method applies to the more general case of Curtis-Tits amalgams.
\ere

\section{{Construction of the amalgam \texorpdfstring{$\amG^\classmap(\fk)$}{amG class map} inside a {group of Kac-Moody type}}}\label{sec:construction}
In this section we shall prove Corollary~\ref{mth:non-trivial completions}.
\subsection{{Orientable Curtis-Tits amalgams}}\label{subsec:orientable}

\paragraph{Proof of Theorem~\ref{thm:concrete completions}.1}
The second statement follows from Theorem~\ref{thm:KM center}. The main statement follows  directly from Corollary 1.2 of~\cite{BloHof2014b}, which we quote here.
\bco \label{mainthm:OCT=Moufang foundation}
Let $\Gamma$ be a {connected} simply laced {Dynkin} diagram with no triangles and $\fk$ a field with at least 4 elements.
The universal completion of {the standard}  Curtis-Tits amalgam over a commutative field $\fk$ and diagram $\Gamma$ is a {group of Kac-Moody type} over $\fk$ with {Dynkin} diagram $\Gamma$ (and $\amG$ is the Curtis-Tits amalgam for this group) if and only if $\amG$ is orientable.
\eco
Recall from Definition~\ref{dfn:universal completion of standard amalgam} that the universal completion of this standard amalgam $\amG$ of type $\Gamma_\classmap(\fk)$ is called the {simply-connected Curtis-Tits group of type $\Gamma_\classmap(\fk)$}. }
In order to fix notation for future reference we recall some relevant details from its proof here.
Fix $\Gm$, $I$, $\fk$, and $\classmap$, as well as $\amG$ as in the Classification Theorem and see Definition~\ref{dfn:standard CT structure} for a precise definition of a standard Curtist-Tits amalgam.

For each $i\in I$, using the identification $\amgrpG_i\cong \SL_2(\fk)$, for each $t\in \fk$ define 
\begin{align}\label{eqn: phi_alpha}
\amgrpU_i^+(t)=\begin{pmatrix} 1 & t \\ 0 & 1 \end{pmatrix} \mbox{ and }
\amgrpU_i^-(t)=\begin{pmatrix} 1 & 0\\ t & 1 \end{pmatrix}%\nonumber
\end{align}
and let \nom{$\amgrpU_i^\vep=\{ \amgrpU_i^\vep(t)\colon t\in \fk^+\}\le\amgrpG_i$}{}, for $\vep=+,-$.
Since $\im \classmap\le \Aut(\fk)$ (i.e.~$\amG$ is orientable) we have that $\amg_{i,j}(\amgrpU_i^+)$ and $\amg_{j,i}(\amgrpU_j^+)$ belong to the same Borel subgroup of $\amgrpG_{i,j}$ for all $i,j\in I$.
This allows us to create a sound Moufang Foundation whose rank $1$ and rank $2$ residues are buildings associated to split algebraic groups.  Then we use the integration results by M\"uhlherr~\cite{Mu1999} to prove that there exists a Moufang twin-building
\nom{$\De=(\De_+,\De_-,\delta_*)$}{} with diagram $\Gm$ with this foundation. This building is constructed from a combinatorial building over the covering of $\Gamma$ associated to $\pi_1(\Gamma,i_0)$ by taking certain fixed points under the action of $\pi_1(\Gamma,i_0)$. Note that M\"uhlherr's theorem allows infinite type sets and as such the image of $\omega$ does not have to be finite. This building has a pair of opposite chambers \nom{$c_+\in \De_+$}{}, \nom{$c_-\in \De_-$}{} such that $\amgrpU_i^\vep$ can be identified with the root group \nom{$\compU_{\vep\alpha_i}\le \Aut(\De)$}{} ($\vep=+,-$, $i\in I(\Gamma)$). Here \nom{$\Pi=\{\alpha_i\colon i\in I\}$}{} is a fundamental system of a system \nom{$\Phi$}{} of (real) roots of type $\Gamma$.
That is, if  \nom{$\Sigma=(\Sigma_+,\Sigma_-)$}{} is the twin-apartment of $\De$ determined by $(c_+,c_-)$, then 
 $\Pi$ is the set of roots in $\Sigma_+$ containing $c_+$ and determined by the panels on $c_+$.
Moreover, let \nom{$\KM=\langle \compU_\alpha\colon \alpha\in \Phi\rangle\le \Aut(\De)$ }{} (In fact since $\Gamma$ is simply-laced we have 
 $\KM=\langle \compU_\alpha, \compU_{-\alpha}\colon \alpha\in \Pi\rangle$).
Then, it is known that $(\KM,\{\compU_\alpha\}_{\alpha\in \Phi})$ is a group with twin root group datum. 
For $\alpha,\beta\in \Pi$, let 
\begin{align*}
\KM_\alpha&=\langle \compU_\alpha,\compU_{-\alpha}\rangle,\\
\KM_{\alpha,\beta}&=\langle \KM_\alpha,\KM_\beta\rangle.
\end{align*}
\idx{$\KM_\alpha$}\idx{$\KM_{\alpha,\beta}$}
The proof of Corollary~\ref{mainthm:OCT=Moufang foundation} furthermore asserts that the identification isomorphisms   
 $\amgrpU_i^\vep\to \compU_{\vep\alpha_i}$ extend to surjective  homomorphisms
  \begin{align*}
\compg_i\colon \amgrpG_i & \stackrel{\cong}{\to} \KM_{\alpha_i}\\
\compg_{i,j}\colon \amgrpG_{i,j} & \to \KM_{\alpha_i,\alpha_j} 
 \end{align*}
where the latter is an isomorphism whenever $\{i,j\}$ is an edge of $\Gamma$, provided one assumes that $|I|\ge 4$.
Universality of $(\ucompG,\ucompg_\bullet)$ gives that the collection \nom{$\compg_\bullet=\{\compg_i,\compg_{i,j}\colon i,j\in I\}$}{} induces a surjective homomorphism $\ucompG\to \KM$. Thus, $(\KM,\compg_\bullet)$ is a completion of $\amG$.
Combined with Theorem~\ref{thm:KM center}, this completes the proof of part 1 in Theorem~\ref{thm:concrete completions}.\epf

\medskip

For future reference also recall~\cite{Re2004} that, setting 
\begin{align*}
\compU^\vep&=\langle \compU_\alpha\colon \alpha\in \Phi^\vep\rangle &&  (\vep=+,-) \\
\compN&=\langle m(u)\colon u\in \compU_\alpha-\{1\}, \alpha\in \Pi\rangle,\\
\compD&=\bigcap_{\alpha\in \Phi}N_{\KM}(\compU_\alpha)
\end{align*}
 we have a twin- BN-pair $((\compB^+,\compN), (\compB^-,\compN))$, where $\compB^\vep=\compD\compU^\vep$. In fact this is the twin BN-pair giving rise to $\De=(\De_+,\De_-,\delta_*)$.
\idx{$\compU^\vep$}\idx{$\compN$}\idx{$\compD$}\idx{$\compB^\vep$}
\subsection{{Non-orientable Curtis-Tits amalgams and the proof of Theorem~\ref{thm:concrete completions}.2.}}\label{subsec:non-orientable CT groups}\label{sec:twisted amalgam}
From now on we fix non-orientable $\delta$ and drop it as a superscript.

\paragraph{Covering graphs, fibers and deck-transformations}
Let \nom{$K$}{} be the kernel of the composition of homomorphisms $\Pig\stackrel{\classmap}{\to} \Aut(\fk)\times\langle \tau\rangle\to \langle \tau\rangle$.
Recall that there is a covering graph 
\nom{$\La=(\hI,\hE)$}{} such that if \nom{$p\colon \La\to \Gamma$}{} is the projection map and $\hat{0}$ is a vertex mapping to $0$, the induced homomorphism is an isomorphism \nom{$p_*\colon \pi_1(\La,\hat{0})\to K$}{} and as $K$ is normal, \nom{$\Pig/K$}{} acts as a group of deck transformations commuting with $p$.
Let \nom{$\theta$}{} be the generator of $\Pig/K\cong \ZZ_2$.

Let $i,j\in I$ be distinct, and let \nom{$\La_{i,j}$}{} denote the subgraph induced on the $p$ fiber over $\{i,j\}$.
Then, since $\Gamma$ has no circuits of length $\le 2$,  $\La_{i,j}$ consists of pairwise disjoint edges if $\{i,j\}\in E(\Gamma)$ and has no edges otherwise.
Since $K\normal \Pig$, in fact $\theta$ restricts to an automorphism of $\La_{i,j}$.

\bde\label{dfn:hamG}
Let \nom{$\amL=\{\amgrpL_{i},\amgrpL_{i, j}, \aml_{i,j} \mid  i, j \in \hI\}$}{} be the amalgam such that, for all $i,j\in \hI$, 

\begin{enumerate}
\item[($\amL$1)] 
$\amgrpL_i\mbox{  is a copy of }\amgrpG_{p(i)}$,
\item[($\amL$2)] 
$ \amgrpL_{i,j}\mbox{  is a copy of }
\begin{cases}
\amgrpG_{p(i),p(j)} &\mbox{if  } \{i,j\}\in \hE, \\
\amgrpG_{p(i)}\times \amgrpG_{p(j)} & \mbox{else. } \\
\end{cases}$
\item[($\amL$3)] 
$ \aml_{i,j}
 =
 \begin{cases}
\amg_{p(i),p(j)} & \mbox{if } \{i,j\}\in \hE,\\
 \mbox{ canonical inclusion } & \mbox{else.}
\end{cases}$
\end{enumerate}

Given any $J\sbe \hI$ with $1\le |J|\le 2$, and denoting by \nom{$\ammapp_\bullet\colon \amL\to \amG$}{} the homomorphism of amalgams induced by $p$, then we have a commuting diagram of isomorphisms

\begin{align}\label{eqn:hamG covers amG}
\xymatrix{
x\in \amgrpL_J\ar[dr]^{\ammapp_J} \ar[rr]^{\theta_J} && \amgrpL_{\theta(J)} \ni x \ar[dl]^{\ammapp_{\theta(J)}}\\
                            & x\in \amgrpG_{p(J)} & 
}
\end{align}
That is, identifying $\amgrpG_{p(J)}$ with its copies $\amgrpL_J$ and $\amgrpL_{\theta(J)}$, the maps \nom{$\ammapp_\bullet$}{} and \nom{$\theta_\bullet$}{} are given by the identity mapping. Then, by ($\amL$3) $\theta_\bullet$ is an automorphism of $\amL$.
\ede
Let \nom{$\widehat{\classmap}=\classmap\after p_*$}.

\bde\label{dfn:hat RGD}
We first apply the construction of Subsection~\ref{subsec:orientable} replacing $\Gamma$, 
 $\classmap$, $\amG$, $\compG_\bullet$, and $(\compG,\compg_\bullet)$ by $\La$, $\widehat{\classmap}$, $\amL$, $\compL_\bullet^\rsc$, and $(\compL^\rsc,\compl_\bullet^\rsc)$ respectively, and leaving all other notation as it is.
Thus $\compL^\rsc$ is the simply-connected Curtis-Tits group of type $\La_{\widehat{\classmap}}(\fk)$.
Next, we will pass to the adjoint setting via the canonical map $\compL^\rsc\to\compL^\rsc/Z(\compL^\rsc)$ replacing 
 $\compL^\rsc_\bullet$ and $(\compL^\rsc,\compl_\bullet^\rsc)$ by $\compL_\bullet$ and $(\compL,\compl_\bullet)$, not changing the notation otherwise (so for instance $\compU_\alpha$ now denotes the (isomorphic) image of the corresponding root group in $\compL^\rsc$). 
Thus, \nom{$(\compL,\{\compU_\alpha\}_{\alpha\in \Phi})$}{} is a twin-root group datum, where \nom{$\Phi$}{} is the root system of type $\La$ with
 fundamental system \nom{$\hPi=\{\alpha_i\mid i\in \hI\}$}{} and 
for each $\alpha,\beta\in \Phi$ we have $\compL_\alpha=\langle \compU_\alpha, \compU_{-\alpha}\rangle$ and $\compL_{\alpha,\beta}=\langle \compL_\alpha, \compL_\beta\rangle$. 
It follows that, not only $(\compL^\rsc,\compl_\bullet^\rsc)$, but also $(\compL,\compl_\bullet)$ is a non-trivial completion of $\amL$.
Now let 
\begin{align*}
\ocompD_\alpha&=N_{\compL_\alpha}(\compU_\alpha)\cap N_{\compL_\alpha}(\compU_{-\alpha}),\\
\ocompD&=\prod_{i\in \hI} D_{\alpha_i}.
\end{align*}
\ede
The group $\compL$ acts on the twin-building $\De=(\De_+,\De_-,\delta_*)$ associated to the twin-root datum, alternatively described as in~\cite[\S 5]{BloHof2014b}.
The automorphism $\theta_\bullet$ induces a diagram automorphism of the root group data of $\compL$.
Indeed,  if \nom{$(\ucompL,\ucompl)$}{} is the universal completion of $\amL$, then $\theta_\bullet$ induces an automorphism of $\ucompL$. Hence $\theta_\bullet$ induces an automorphism of $\compL\cong \ucompL/Z(\ucompL)$. 
More precisely, the following diagrams commute
 \begin{align}\label{eqn:phi theta=theta phi}
 \xymatrix{
 \amgrpL_i \ar[r]^{\compl
_i}\ar[d]_{\theta_i} & \compL\ar[d]^{\theta} \\
 \amgrpL_{\theta(i)} \ar[r]^{\compl
_{\theta(i)}} & \compL\\
}
&&
\xymatrix{
 \amgrpL_{i,j} \ar[rr]^{\compl
_{i,j}}\ar[d]_{\theta_{i,j}} && \compL\ar[d]^{\theta} \\
 \amgrpL_{\theta(i),\theta(j)} \ar[rr]^{\compl
_{\theta(i),\theta(j)}} && \compL.\\
}
\end{align}
Note that $\theta(\compU_{\alpha_i})\in \{\compU_{\alpha_i}, \compU_{-\alpha_i}\}$. Moreover, if $\{i,j\}\in \hE$, then because of the Chevalley relations, if  $\theta(\compU_{\alpha_i})=\compU_{\epsilon\alpha_i}$, then also $\theta(\compU_{\alpha_j})=\compU_{\epsilon \alpha_j}$, for some $\epsilon=\pm$.
Since $\amG$ is a non-orientable amalgam, there exists some $i$ for which $\theta(\compU_{\alpha_i})=\compU_{-\alpha_i}$, and so by connectivity of $\La$, we must have $\theta(\compU_{\alpha_i})=\compU_{-\alpha_i}$, for all $i\in \hI$. This shows that $\theta$ interchanges the two halves of $\De$ and there exists a chamber $c_+\in \De_+$, such that $c_+\opp c_+^\theta$.
Namely, $c_+$ is the chamber corresponding to the Borel subgroup corresponding to the fundamental system $\Pi$.

\medskip
Let \nom{$\compLth$}{} denote the group of fixed points of $\compL$ under $\theta$. Our next aim is to show that $\compLth$ contains a non-trivial completion of $\amG$.

\medskip

We consider the following setup.
Let $i,j\in I(\Gamma)$ and $k\in p^{-1}(i)$ and $l\in p^{-1}(j)$ so that if $\{i,j\}\in E(\Gamma)$, then $\{k,l\}\in E(\La_{i,j})$.

\ble\label{lem:theta orbit groups}
\hspace{2em}
\begin{enumerate}
\item 
$\compl
(\amgrpL_{k,\theta(k)})\cong \compl
(\amgrpL_k)\times \compl
(\amgrpL_{\theta(k)})\cong \SL_2(\fk)\times\SL_2(\fk)$,
\item
$\langle \compl
(\amgrpL_{k,l}),\compl
(\amgrpL_{\theta(k),\theta(l)})\rangle \cong \compl
(\amgrpL_{k,l})\circ \compl
(\amgrpL_{\theta(k),\theta(l)})$.
\end{enumerate}
\ele
\bpf
By Lemma~\ref{lem:vertex and edge group images}  $\compl
(\amgrpL_k)\cong\SL_2(\fk)\cong\compl
(\amgrpL_{\theta(k)})$.
Moreover, since $k$ and $\theta(k)$ are at distance more than $2$ in $\La$, and $\Gamma$ is triangle free and connected,  by Lemma~\ref{lem:A1 x A2},  $\compl
(\amgrpL_{k,\theta(k)})\cong \compl
(\amgrpL_k)\times \compl
(\amgrpL_{\theta(k)})$.

Similarly, it follows from the structure of $\La_{i,j}$ that $\compl
(\amgrpL_{k,l})$ and $\compl
(\amgrpL_{\theta(k),\theta(l)})$ commute, so they can only intersect in their center. 
\epf
\medskip

\bde\label{dfn:hamGtheta}
Next, define the amalgam 
\nom{$\amL^{\theta}=\{\amgrpL_{i}^{\theta},\amgrpL^{\theta}_{i, j}, \aml^{\theta}_{i,j} \mid  i, j \in I\}$}{}, where
for $i,j\in I$ and $k,l\in \hI$ as before, we define the following subgroups of $\compL$:
\begin{align*}
\amgrpL^\theta_i = \{\compl
_k(y)\theta(\compl
_{k}(y)) \mid y\in \amgrpL_{k}\} &&\mbox{ and }&&
 \amgrpL^\theta_{i,j} =\{\compl
_{k,l}(y)\theta(\compl
_{k,l}(y)) \mid y\in \amgrpL_{k,l}\},
\end{align*}
where the connecting maps $\aml_{i,j}$ are inclusions of subgroups of $\compL$.
Also define a map of amalgams \nom{$\ammapm_\bullet\colon \amG\to \amL^\theta$}{} setting
\begin{align*}
\ammapm_i\colon y & \mapsto \compl
_k(y)\theta(\compl
_{k}(y)) &&\mbox{ for all }y\in \amgrpG_i,\\
\ammapm_{i,j} \colon y & \mapsto \compl
_{k,l}(y)\theta(\compl
_{k,l}(y))  &&\mbox{ for all }y\in \amgrpG_{i,j}.
\end{align*}
Note that we are implicitly using the map $\ammapp_\bullet$.
 \ede
\bpr\label{prop:amG delta in amG delta K}
\hspace{2em}
\begin{enumerate}
\item
 $\amL^\theta$ is an amalgam of subgroups in $\compLth$.
\item $\ammapm_\bullet\colon\amG\to \amL^\theta$ is a non-trivial morphism of amalgams.\end{enumerate}
\epr
\bpf
Part 1. In the notation of Definition~\ref{dfn:hamGtheta}, for any $y\in \amgrpG_J$, we have $\theta\big(\compl_J(y)\theta(\compl_{J}(y))\big)
=\theta(\compl_J(y))\compl_{J}(y)=\compl_J(y)\theta(\compl_{J}(y))$ for 
 $J=\{k\}$ and $J=\{k,l\}$ due to the central and direct products appearing in Lemma~\ref{lem:theta orbit groups}.

Part 2. 
To see that $\ammapm_i$ is non-trivial, we show that $\ker\ammapm_i$ is central in $\amgrpG_i$.
Namely, if $y\in \amgrpG_i$ is such that $\ammapm_i(y)=1$, then apparently $\compl_k(y)\in \amgrpL_k\cap \amgrpL_{\theta(k)}$ which is central in $\amgrpL_k$ (and $\amgrpL_{\theta(k)}$).
Since the map $y\mapsto y$ identifies $\amgrpG_i$ with $\amgrpL_k$ the claim follows.

Since the products appearing in Lemma~\ref{lem:theta orbit groups} are direct or central, $\ammapm_i$ and $\ammapm_{i,j}$ are group homomorphisms.
Hence it suffices to show that
$\ammapm_{i,j}\after \amg_{i,j}=\aml^\theta_{i,j}\after\ammapm_i$.
To this end note that the connecting maps of $\amL$ are compatible with the identification in~\eqref{eqn:hamG covers amG} of $\amgrpL_k$ with $\amgrpG_i$ (resp. $\amgrpL_{k,l}$ with $\amgrpG_{i,j}$) and  of Equation~\eqref{eqn:phi theta=theta phi}.
\epf

\paragraph{Proof of Theorem~\ref{thm:concrete completions}.2 }
By Proposition~\ref{prop:amG delta in amG delta K}, the non-orientable amalgam $\amG$ has a non-trivial completion in the group $\compL^\theta$. This proves part 2.~of Theorem~\ref{thm:concrete completions}.

\medskip
From now on we shall denote by \nom{$(\compG,\compg_\bullet)$}{} the completion of $\amG$ generated by $\amL^\theta$ in $\compLth$.
\medskip

\subsection{Proof of Theorems~\ref{mth:non-trivial completion iff sim_A in B} and Theorem~\ref{mth:amA injects iff sim_1 refines sim_A}}\label{subsec:non-trivial completion iff sim_A in B}
\paragraph{Proof of Theorem~\ref{mth:non-trivial completion iff sim_A in B}}
Let $\amA$ be any Curtis-Tits amalgam of type $\Gamma_\classmap(\fk)$ and let $(\ucompA,\ucompa_\bullet)$ be its universal completion.
Since $\amG_\classmap(\fk)$ is a universal amalgam, we have a map $\rho_\bullet\colon \amG_\classmap(\fk)\stackrel{\uammapp_\bullet}{\to} \amA\stackrel{\ucompa_\bullet}{\to} \ucompA$ and $(\ucompA,\rho_\bullet)$ is a completion of $\amG_\classmap(\fk)$.  
By Proposition~\ref{prop:sim_R refines sim_0}, $\sim_\amA$ refines $\sim_0$.

For the converse, let $\amG_\classmap(\fk)$  be the standard Curtis-Tits amalgam of type $\Gamma_\classmap(\fk)$. 
Then, by Corollary~\ref{mainthm:OCT=Moufang foundation} and Proposition~\ref{prop:amG delta in amG delta K} it has a non-trivial universal completion $(\ucompG,\ucompg_\bullet)$.
Consider the completion $(\compG_0,\compg_{0;\bullet})$, where $\compG_0=\ucompG/Z(\ucompG)$ and its relation 
 $\sim_{\compG_0}$.
 Then, by the preceding paragraph $\sim_{\compG_0}$ refines $\sim_0$.
 We now claim that we have equality. Namely, let $i,j\in I$ be such that $ i\sim_0 j$ and let 
  $z$ be the central involution of $\amgrpG_{i,j}=\amgrpG_i\times \amgrpG_j$.
 Then, by Definition~\ref{dfn:relation 0}, for any $k\in I-\{i,j\}$, either $k$ is not connected to either $i$ or $j$ in $\Gamma$ or it is connected to both. In the former case, clearly $\compg_{0;i,j}(z)$ commutes with $\compg_{0;k}(\amgrpG_k)$ and in the latter, it follows from Lemma~\ref{lem:A3} that $\compg_{0;i,j}(z)$ commutes with $\compg_{0;k}(\amgrpG_k)$. That is, apparently $\compg_{0;i,j}(z)\in Z(\compG_0)=\{1\}$ and we must have $i\sim_{\compG_0} j$ already.
  
It now follows that if $\sim_\amA$ refines $\sim_0$, then we have a homomorphism of amalgams 
 $\amA\to \amG_0$, where $\amG_0$ is the image of $\amG_\classmap(\fk)$ under $\compg_{0;\bullet}$ in $\compG_0$. In particular, $\amA$ has  $\compG_0$ as a non-trivial completion and hence $\amA$ also has a non-trivial universal completion.
 \epf

\paragraph{Proof of Theorem~\ref{mth:amA injects iff sim_1 refines sim_A}}
Suppose $\amA$ is a Curtis-Tits amalgam with property (D) of type $\Gamma_\classmap(\fk)$ such that $\sim_\amA$ refines $\sim_0$.
By Theorem~\ref{mth:non-trivial completion iff sim_A in B} this is equivalent to the assumption that $\amA$ has a non-trivial universal completion, which we shall denote $(\ucompA,\ucompa_\bullet)$.
Since $\amG_{\classmap}(\fk)$ is universal, there is a map $\amG_{\classmap}(\fk)\to \amA$. By universality of the completion $(\ucompG,\ucompg_\bullet)$ of $\amG_{\classmap}(\fk)$, we have a surjective homomorphism $\pi\colon \ucompG\to \ucompA$ taking  $\ucompg_\bullet(\amG_{\classmap}(\fk))=\amG_1$ to $\bar{\amA}$.
Hence, we have
 $d^\ucompA=\pi\after d^\ucompG$. 
 Define $\kIZ_\ucompA=\langle \kIZ^1, \kIz_{kl}\colon k\sim_\amA l\mbox{ with } k,l\in I\mbox{ distinct}\rangle\le (\fk^*)^I$.
 Note that for any $i,j\in I$, 
 \begin{align*}
 i\sim_{\bar{\amA}} j & \iff \kIz_{ij}\in \kIZ^\ucompA,\\
 i\bar{\sim_{\amA}} j & \iff \kIz_{ij}\in \kIZ_\ucompA.
 \end{align*}
Note that we have $\kIZ_\ucompA\le \kIZ^\ucompA$. Consequently, we have 
 $d^{\ucompG}(\kIZ_\ucompA)\le d^{\ucompG}(\kIZ^\ucompA)\le \ker \pi$.
Consider the quotient $\underline{\ucompA}=\ucompG/d^{\ucompG}(\kIZ_\ucompA)$.
Then, we have a canonical map $\underline{\pi}\colon \underline{\ucompA}=\ucompG/d^{\ucompG}(\kIZ_\ucompA)\to \ucompG/d^{\ucompG}(\kIZ^\ucompA)\to \ucompG/\ker \pi=\ucompA$.

Since $\kIz_{ij}\in \kIZ_{\ucompA}$ for all distinct $i,j\in I$ with $i\sim_\amA j$, there is a completion $(\underline{\ucompA},\underline{\ucompa})$ of $\amA$.
Note that $\underline{\pi}$  takes $\underline{\ucompa}(\amA)$ to  $\bar{\amA}$.
By universality of $(\ucompA,\ucompa_\bullet)$ we now must have $\underline{\ucompA}=\ucompA$ and 
 we obtain $\kIZ_\ucompA=\kIZ^\ucompA$ and $\bar{\sim_\amA}=\sim_{\bar{\amA}}$.
The second claim of Theorem~\ref{mth:amA injects iff sim_1 refines sim_A} follows immediately from the first.
\epf

\section{{The geometry \texorpdfstring{$(\De^\theta,\approx)$}{De theta approx} for  non-orientable Curtis-Tits amalgams}}\label{section:De theta}
\subsection{Definition of \texorpdfstring{$(\De^\theta,\approx)$}{De theta approx}}\label{subsec:Detheta definition}
The construction of $\De^\theta$, and the proof of its simple-connectedness follows the pattern of~\cite{BloHof2014a}.

Let $\De=(\De_+,\De_-,\delta_*)$ be the twin-building associated to $\compL$ as in Section~\ref{sec:twisted amalgam} and let $\opp$ denote opposition.
Note that $\De$ has diagram $\La$ over the index set $\hI$.

Recall that the action of $\theta$ on $\compL$ induces an action of $\theta$ on $\De$ that interchanges $\De_+$ and $\De_-$, and such that 
 $c_+\opp c_+^\theta=c_-$, where, for $\vep=+,-$,  $c_\vep$ is the chamber corresponding to the Borel group $B_\vep$, containing $\compU_{\vep\alpha}$ for all $\alpha\in \Pi$.
 
\bde
We define a relaxed incidence relation on $\De_\vep$ as follows.
We say that $d_\vep$ and $e_\vep$ are $i$-adjacent if and only if 
$d_\vep$ and $e_\vep$ are in a common $p^{-1}(i)$-residue.
In this case we write 
$$d_\vep\approx_i e_\vep,$$
where we let $i\in I$.
Note that the residues in this chamber system are $J$-residues of $\De_\vep$ where $J^\theta=J\sbe \hI$.

Define the fixed subgeometry
 $\De^\theta=\{d\in \De_+ \mid d \opp d^\theta\}$ endowed with the adjacency relation $\approx$.
Note that $c_+\in \De^\theta$.

It is easy to see that residues of $\De^\theta$ are the intersections of certain residues of $(\De_+,\approx)$ with the set $\De^\theta$.
\ede 
\bre
The motivation for defining $\De^\theta$ is that the action of  $\compL$ on $\De$ induces an action of  $\compLth$ on $\De^\theta$.
We shall prove in Section~\ref{sec:L=Gtheta} that this action is flag-transitive and use Tits' Lemma~\cite{Ti1986b} to show that under certain conditions, the universal completion of $\amG$ is a central extension of $\compLth$.

\ere

The following results from~\cite{BloHof2014a} and their proofs remain valid for general {connected } simply-laced diagrams without triangles and we mention them here for their value in this more general context:
 Lemma 4.24 (using that $\theta$ has no fixed vertices or edges on $\Lambda$), 4.26, 4.27, Lemma 4.28, Lemma 4.30, Proposition 4.31, Proposition 4.32 and Theorem 4.45. 
As in loc.~cit.~we obtain the following theorem.

\bth\label{thm:Detheta simply connected}
Let $\fk$ be a field of size at least $7$.
Then, $\De^\theta$ is connected and simply $2$-connected.
\eth

The following corollary is of general interest.
{\bco
Let $d\in \De_\vep$.
Then, there exists $(c,c^\theta)\in \De^\theta$ such that 
 $d\in \Sigma_+(c,c^\theta)$.
\eco
\pf
Let $u=\delta_*(d,d^\theta)$. We induct on $l(u)$.
If $u=1$, we are done.
Assume $l(u)>0$. From Lemma 4.26 in ~\cite{BloHof2014a} 
we know that $u=w(w^{-1})^\theta$ is irreducible for some irreducible $w\in W$.
Pick any $s_i\in S$ such that $l(s_iw)=l(w)-1$ and let $e$ be a chamber $i$-adjacent to $d$.
Then, $\delta_*(e,e^\theta)=s_iw(w^{-1})^\theta s_{\theta(i)}$ is shorter by $2$, so by induction there exists
 $c$ such that $e\in \Sigma(c,c^\theta)_+$.
Moreover, calling $\pi$ the $i$-panel on $d$ we have
 $\proj^*_\pi(e^\theta)=d\in \Sigma(c,c^\theta)_+$.
\qed
}
\section{The action of the universal completion \texorpdfstring{$\ucompG$}{tG} on \texorpdfstring{$\De^\theta$}{De theta}}\label{section:action on Detheta}
\subsection{The torus of \texorpdfstring{$\compLth$}{Gth}}\label{subsec:torus}
\bde
Let $\compL$ and $\compLth$, and $\amL^\theta$ be as in Subsection~\ref{sec:twisted amalgam} and let $(\ucompG,\ucompg)$ be the universal completion of  $\amG$ and recall from Proposition~\ref{prop:amG delta in amG delta K}
 that we have a non-trivial morphism of amalgams $\ammapm_\bullet\colon\amG\to \amL^\theta$.
By the universal property we have a map 
\begin{align*}
\umap \colon \ucompG &\rightarrow \compLth
\end{align*}
 that restricts to a non-trivial homomorphism of amalgams on $\amL^\theta$.
Recall that  $\compG=\umap(\ucompG)$ by definition. For $i\in I$, let 
 \nom{$\ucompD_i=\ucompg_i(\amgrpD_i)\le \ucompG$}{}, 
  \nom{$\ucompD=\langle \ucompD_i\colon i\in I(\Gamma)\rangle\le \ucompG$}{}, and 
 \nom{$\compD_i^\theta=\ammapm_i(\amgrpD_i)=\umap(\ucompD_i)$}{}. 
Finally, let \nom{$\compD_{\theta}=\prod_{i\in I} \compD_i^\theta\le \ocompD\cap \compLth$}{} (where $\ocompD$ is as in Definition~\ref{dfn:hat RGD}).
\ede

\bre
All arguments up until Theorem~\ref{thm:kerphi = ZtGcaptD} do not require the assumption that $Z(\compL)=1$ and are valid for any quotient of $\compL^\rsc$ that admits $\theta$ as an automorphism.
\ere

The aim of this subsection is to investigate the difference between $\compD_{\theta}$ and $\ocompD\cap\compLth$.

\medskip
\noindent Lemma 3.3 of~\cite{Cap2007} has the following consequence.
\ble\label{lem:stabilizer of opposite chambers}
We have $\Stab_{\compL}((c_+,c_-))=\bigcap_{\alpha\in\Phi} N_{\compL}(\compU_\alpha)=\prod_{i\in \hI} \ocompD_{\alpha_i}=\ocompD$.
Therefore, $\Stab_{\compLth}((c_+,c_-))=\ocompD\cap\compLth$
 and $\Stab_\compG((c_+,c_-))=\ocompD\cap\compG$.
\ele

\bde
We now define the following maps. 
$d\colon  (\fk^*)^{\hI} \to  \ocompD$, and $\theta, \nu\colon  (\fk^*)^{\hI}\to(\fk^*)^{\hI}$ by 
\begin{align*}
d(x_i)_{i\in \hI}=\prod_{i\in \hI}d_i(x_i); &&
{\theta\left( (x_i)_{i\in \hI} \right)  }= (x_{\theta(i)})_{i\in \hI}; &&
\nu(x)=x \theta(x)^{-1};
\end{align*}
Note that $\theta$ and $\nu$ on $(\fk^*)^{\hI}$ induce $\theta$  and $\nu$ on $\ocompD$ via $d$ (see~\eqref{eqn:phi theta=theta phi}).
We now set 
 \begin{align*}
\kIF&=\ker \nu=\{x\in (\fk^*)^{\hI} \mid x=\theta(x)\}, \\
 \imnu&=\im \nu, \\
 \kIK&=\ker d.
 \end{align*}
\ede

Now note that $d(\kIF)=\compD_{\theta}=\langle d_i(a)d_{\theta(i)}(a)\mid i\in \hI, a\in \fk^*\rangle\le \ocompD\cap \compLth$.
However, this inclusion may be proper, as we will now investigate.
The diagram in Figure~\ref{fig:snake} arises from applying the snake lemma to the commuting diagram involving the two rows connected via $d$, which are exact.
The columns are all exact.

\begin{figure}
\begin{centering}
\begin{tikzpicture}
  \matrix[matrix of math nodes,column sep={50pt,between origins},row
    sep={45pt,between origins},nodes={rectangle}] (s)
  {
      &|[name=zFK]| 0 &|[name=zK]| 0 &|[name=zNK]| 0  \\
|[name=hzFK]| 0       &|[name=FK]| \kIF\cap \kIK &|[name=K]| \kIK &|[name=NK]| \imnu\cap \kIK & |[name=hzNK]|     \\
  |[name=hzF]| 0      &|[name=F]| \kIF &|[name=k2n]| (\fk^*)^{2n} &|[name=N]| \imnu &|[name=hzN]| 0 \\
    |[name=hzDGtheta]| 0 &|[name=DGtheta]| \ocompD\cap \compLth &|[name=D]| D &|[name=imnu]| \nu(\ocompD) &|[name=hzimnu]| 0 \\
    &|[name=zDGtheta]| (\ocompD\cap \compLth)/\compD_{\theta} &|[name=zD]| 0 &|[name=zimnu]| 0 \\
  };
  \draw[->] (zFK) edge (FK)
            (zK) edge (K)
            (zNK) edge (NK)
            (hzFK) edge (FK)
            (FK) edge (K)
            (K) edge node[auto] {$\nu$} (NK)
            (FK) edge  (F)
            (K) edge  (k2n)
            (NK) edge  (N)
            (hzF) edge (F)
            (F) edge (k2n)
            (k2n)  edge node[auto] {$\nu$} (N)
            (N)  edge (hzN)
             (F) edge node[auto] {$d$} (DGtheta)
            (k2n)  edge node[auto] {$d$} (D)
            (N)  edge node[auto] {$d$} (imnu)
            (hzDGtheta) edge (DGtheta)
            (DGtheta)  edge (D)
            (D)  edge node[auto] {$\nu$} (imnu)
            (imnu) edge (hzimnu)
            (DGtheta)  edge (zDGtheta)
            (D)  edge (zD)
             (zDGtheta) edge (zD)
            (zD) edge (zimnu)
        (imnu) edge (zimnu)
  ;
  \draw[->,gray,rounded corners] (NK) -| node[auto,text=black,pos=.7]
    {\(\partial\)} ($(hzN.east)+(.5,0)$) |- ($(k2n)!.35!(D)$) -|
    ($(hzDGtheta.west)+(-.5,0)$) |- (zDGtheta);
\end{tikzpicture}
\caption{The $\ocompD\cap \compLth$-$\compD_{\theta}$ diagram.}\label{fig:snake}
\end{centering}
\end{figure}
The snake lemma gives the following result.
\ble\label{lem:snake}
We have an isomorphism
 $\partial\colon (\imnu\cap \kIK)/\nu(\kIK)\cong (\ocompD\cap \compLth)/\compD_{\theta}$.
\ele
We shall show in Corollary~\ref{cor:Gtheta=L Stabc} that $\compG=\compLth$ if and only if the following condition holds:
\begin{align}\tag{D}
(\ocompD\cap \compLth)/\compD_{\theta}=0  \label{eqn:DGtheta=Dtheta}
\end{align}
By Lemma~\ref{lem:snake}, this is equivalent to the fact that $\kIM\cap\kIK=\nu(\kIK)$.

\ble\label{lem:DGtheta mod Dtheta}
{ $(\imnu\cap \kIK)/ \nu(\kIK)\cong (\ocompD\cap \compLth)/\compD_{\theta}$ is an elementary abelian $2$-group.
 If $\fk\cong \FF_q$, then it has rank at most $n$.}
\ele
\bpf
Let $x=(x_i)_{i\in \hI}\in \imnu\cap \kIK$. This means that 
  $x_i=a_ia_{\theta (i)}^{-1}=\nu(a)$ for some $a=(a_i)_{i\in \hI}$
   and $d(x)=1\in \ocompD$.
Note that $x^2=x \theta(x)^{-1}=\nu(x)\in \nu(\kIK)$.
This means that $(\imnu\cap \kIK)/\nu(\kIK)$ is an elementary abelian $2$-group. 
{It is a quotient of a subgroup of $\imnu$ which, if $\fk=\FF_q$,  is a product of $n$ cyclic groups isomorphic to $\FF_q^*$, and the second claim follows.}
\epf

\bre
If we define $\compL$ so that $Z(\compL)=1$ (see Definition~\ref{dfn:hat RGD}), we have $\kIK=\ker\comop_\amL$.
In this case we can verify condition~\eqref{eqn:DGtheta=Dtheta} combinatorially (see Example~\ref{example verifying condition D}).
\ere

\subsection{The amalgam of parabolics}\label{subsec:parabolic amalgam}
For any subset $J\sbe I$, let $\hJ=p^{-1}(J)\sbe \hI$, where $p\colon \La\to \Gamma$ is the projection map.

\bde
Note that $\compG$ is the subgroup of $\compLth$ generated by the $\umap$-images of the groups in $\amL^\theta$.
Let \nom{$\amB=\{B_J\mid J\sbe I\}$}{} denote the amalgam of all parabolic subgroups of $\compG$ with respect to its  action on $\De^\theta$. 
That is, for any $J\sbe I$, we let 
$B_J=\Stab_\compG(R_\hJ)$, where $R_\hJ$ is the $\hJ$-residue of $\De_+$ on $c_+$ such that $R_\hJ^\theta \opp R_\hJ$.  
In particular, $B_\emptyset=\ocompD\cap\compG$.
Also, for each $J\sbe I$, let \nom{$\compG_J=\langle \umap(\amgrpL^\theta_j)\mid j\in J\rangle_{\compG}$}{} (see Definition~\ref{dfn:hamGtheta}).
\ede

\ble\label{lem:amalgam transitive on panels}
The stabilizer $B_{\{i\}}$ in $\compG$ of the $\{i,\theta(i)\}$-residue on $c_+$ of $\De_+$ is transitive on its intersection with $\De^\theta$.
\ele

\bpf
The $\{i,\theta(i)\}$-residue $R$ containing the chamber $c_+$ with $(c_+,c_+^\theta)\in \De^\theta$ has type $A_1\times A_1$.
Therefore by Lemma~\ref{lem:theta orbit groups} the stabilizer in $\compL$ of $R$  acts as $\amgrpG_i\times \amgrpG_{\theta(i)}\cong \SL_2(\fk)\times \SL_2(\fk)$ on $R$ and the stabilizer in $\compG$ of $R$ contains $\umap(\amgrpL^\theta_i)= \{\umap(gg^\theta)\mid g\in \amgrpL_i\}$.
Let $R_i,R_{\theta(i)}$ be the $i$ and $\theta(i)$-panels on $c^+$.
Identify $R$ with $R_i\times R_{\theta(i)}$ via 
 $d\mapsto (\proj_{R_i}(d),\proj_{R_{\theta(i)}}(d))$.
For any of these chambers $d$ we have
 $\delta_*(d,d^\theta)\in \{1, s_is_{\theta(i)}\}$.

Let $S=\{d\in R\mid \delta_*(d,d^\theta)=s_is_{\theta(i)}\}$.
Note that every panel in $R$ meets $S$ in exactly one chamber.
Also $\umap(\amgrpL^\theta_i)$ acts $2$-transitively on $S$.
Any chamber $c'\in R-S=R\cap \De^\theta$, corresponds to a unique ordered pair of distinct chambers of $S$, that is $(R_i\cap S,R_{\theta(i)}\cap S)$, where $R_i$ is the $i$-panel on $c'$ and $R_{\theta(i)}$ is the $\theta(i)$-panel on $c'$.
By $2$-transitivity $\umap(\amgrpL^\theta_i)$ is transitive on such pairs, hence is transitive on $R\cap \De^\theta$.
\epf

\bco\label{cor:Gtheta=L Stabc}\label{cor:BJ flag transitive}
Suppose that $\De^\theta$ is connected.
For any $J\sbe I$, the group $B_J$ is flag-transitive on 
 $R_{\hJ}\cap \De^\theta$.
In particular, $\compG$ itself is flag-transitive on $\De^\theta$, and as a consequence, so is $\compLth$.
Hence, for any $J\sbe I$, we have $B_J=\compG_J\cdot (\ocompD\cap\compG)$ and $\compLth=\compG \cdot \Stab_{\compLth}((c_+,c_-))=\compG\cdot (\ocompD\cap\compLth)$ (product of groups).
\eco
\bpf
Lemma~\ref{lem:amalgam transitive on panels} combined with the connectedness of $\De^\theta$ in fact shows that already $\compG_J$ is transitive on the chambers of $R_{\hJ}\cap \De^\theta$ and the conclusion follows.
\epf
\medskip

\paragraph{Proof of Theorem~\ref{thm:G flag-transitive on Detheta}}
By assumption of Part 2.~of Theorem~\ref{thm:concrete completions} we have $|\fk|\ge 7$ and so by 
 Theorem~\ref{thm:Detheta simply connected},  $\De^\theta$ is connected and simply $2$-connected. That $\compG$ and hence $\ucompG$ is flag-transitive on $\De^\theta$ follows from Corollary~\ref{cor:BJ flag transitive}.

\bpr\label{prop:universal completion of parabolics}
Suppose that the conclusion of  Theorem~\ref{thm:Detheta simply connected} holds.
For any $J\sbe I$ with $|J|\ge 3$, $B_J$ is the universal completion of the amalgam $\amB_J=\{B_K\mid K\sbne J\}$.
In particular, $\compG$ is the universal completion of the amalgam $\amB$.
\epr
\bpf
This follows from Tits' Lemma~\cite[Corollaire 1]{Ti1986b} combined with Corollary~\ref{cor:BJ flag transitive}.
\epf

\bpr\label{prop:thm B part 3}
$\compG=[\compLth,\compLth]$ and $\compLth/\compG$ is an elementary abelian $2$-group. If $\fk=\FF_q$, then it has rank at most $n=|I|$.
\epr
\bpf
As $\compG$ is generated by perfect groups it is perfect itself, so $\compG\le [\compLth,\compLth]$.
By Corollary~\ref{cor:Gtheta=L Stabc} and Lemma~\ref{lem:stabilizer of opposite chambers} $\compLth = \compG (\ocompD\cap \compLth)$ and the latter group normalizes $\compG$;
namely $\ocompD$ normalizes $\umap(\amgrpL_{i,\theta(i)})$ and since $\compLth\cap \umap(\amgrpL_{i,\theta(i)})= \umap(\amgrpL^\theta_{i})$ we see that $\ocompD\cap \compLth$ normalizes $\umap(\amgrpL^\theta_i)$ for every $i\in I$.
Therefore, $\compLth/\compG$ is isomorphic to a quotient of $(\ocompD\cap \compLth)/\compD_{\theta}$, which is abelian, so that $\compG=[\compLth,\compLth]$. The result follows from Lemma~\ref{lem:DGtheta mod Dtheta}.
\epf

\medskip
\paragraph{Proof of Theorem~\ref{thm:concrete completions}.3}
This follows from Proposition~\ref{prop:thm B part 3} and Theorem~\ref{thm:KM center}.
For a concrete example see Subsection~\ref{example verifying condition D}.
\subsection{\texorpdfstring{$\ucompG$}{tG} is a central extension of \texorpdfstring{$\compLth$}{Gth}}\label{sec:L=Gtheta}
From now on we shall assume that $(\imnu\cap \kIK)={\nu(\kIK)}$. This happens frequently, as $\kIK$ itself is generally already very small.  By Lemma~\ref{lem:snake} and Corollary~\ref{cor:Gtheta=L Stabc} we have 
\begin{align*}
\compG=\compLth \mbox{ and }B_\emptyset=\ocompD\cap\compG=\compD_{\theta}.
\end{align*}
This condition is necessary to ensure that the groups $B_J$ are entirely obtainable from the amalgam $\amL^\theta$ via Corollary~\ref{cor:BJ flag transitive}.

We also assume throughout this section that the conclusion of Theorem~\ref{thm:Detheta simply connected} holds, that is, $(\De^\theta,\approx)$ is simply-connected.
Recall that  $(\ucompG,\ucompg)$ is the universal completion of $\amG$.
We denote proper subgroups of $\ucompG$ and $\compLth$ by $\tH$ and  $H^{\theta}$ respectively, with the assumption that  $\umap(\ucompH)=H^\theta$.
In particular, for $\amgrpG_i\in \amG$ we set $\ucompg(\amgrpG_i)=\ucompG_i$.

\bde
Define the following amalgam in $\ucompG$.
$\tamB=\{\ucompB_J\mid J\sbe I\}$, where
 $\ucompB_J=\langle \ucompG_j, \ucompD_i\mid j\in J, i\in I\rangle_\ucompG$.
The connecting maps of $\tamB$ are given by inclusion of subgroups in $\ucompG$.
Recall that we have a surjective homomorphism 
$\umap\colon \ucompG\to \compG$ given by universality of $\ucompG$. 
For each $J\sbe I$, we also define
 \begin{align*}
 \ucompK_J =\ker\umap\cap \ucompB_J\mbox{ and }
 \ucompK=\langle \langle \ucompK_J\mid J\sbe I\rangle\rangle,
 \end{align*}
 that is, the smallest normal subgroup of $\ucompG$ containing all $\ucompK_J$.
 \ede
\ble\label{lem:tK=ker phi}
For any $J\sbe I$, of rank at least $3$, $\ucompK_J$ is the normal closure in $\ucompB_J$ of all $\ucompK_{J'}$ with $J'\sbne J$.
In particular $\ucompK=\ker\umap$.
\ele

\bpf
By Corollary~\ref{cor:BJ flag transitive}, since $B_J=G_J\cdot (\ocompD\cap\compG)=G_J\cdot \compD_{\theta}$, the restriction  $\umap|_{\ucompB_J}\colon \ucompB_J\to B_J$ takes 
 $\ucompB_{J'}$ onto $B_{J'}$ for any ${J'}\subsetneq J$.
By Proposition~\ref{prop:universal completion of parabolics}, $B_J$ is the universal completion of the amalgam $\{B_{J'}\mid {J'}\sbne J\}$.
Therefore by Lemma~\ref{lem:locally generated kernel},  $\ucompK_J$ is the normal closure in $\ucompB_J$ of all $\ucompK_{J'}$ with $J'\sbne J$.
In particular $\ucompK=\ker\umap$.
\epf

\medskip
Recall that $\ucompD=\langle \ucompg_i(\amgrpD_i)\colon i\in I\rangle$.
Our next aim is to show that $\ucompK\le Z(\ucompG)\cap \ucompD$.
To this end it suffices to show that  $\ker\umap\cap \ucompB_J\le Z(\ucompG)\cap \ucompD$.
We start at rank $0$.
\ble\label{lem:center lemma}\label{lem:rank 0 center lemma}
We have  $\ucompD \cap \umap^{-1}(Z(\compG))\subseteq Z(\ucompG)\cap\ucompD$. 
In particular, $\ucompK_\emptyset=\ucompD\cap \ker \umap\le Z(\ucompG)\cap \ucompD$.

\ele 
\bpf
Note that $\ucompD$ normalises every $\ucompG_i$.
Thus, if $x\in \ucompD\cap \umap^{-1}(Z(\compG))$, then for every $g \in \ucompG_i$, also $g^x \in \ucompG_i$ and $\umap (g) =\umap(g^{x})$. As $\umap |_{\ucompG_i}$ is injective we get $g^x=g$, so $x \in C_{\ucompG}(\ucompG_i)$. This holds for all $i$ so $x\in Z(\ucompG)\cap \ucompD$. \epf

\ble\label{lem:tK rank 2}
For every $J\sbe I$ with $|J|\le 2$ we have 
 $\ucompK_J\le Z(\ucompG)\cap \ucompD$.
\ele
\bpf
Let $x\in \ucompK_J$, then $x=gd^{-1}$ for some $g\in \ucompG_J$ and $d\in \ucompD$.
We now have $\umap(g)=\umap(d)$.
Moreover, for any $y\in \ucompG_J$, we have 
 $\umap(y^g)=\umap(y^d)$, where $y^g,y^d\in \ucompG_J$, as the latter is normalized by both $g$ and $d$.
By Lemma~\ref{lem:vertex and edge group images}, if $J$ is a vertex or an edge, then $\umap|_{\ucompG_J}$ is injective and so we have $y^g=y^d$, so $g\in \ucompG_J$ acts as a diagonal, necessarily inner, automorphism of $\ucompG_J$. Therefore $g\in \ucompD_J$ so that $x\in \ucompD\cap \ker\umap$, which is contained in $Z(\ucompG)$ by Lemma~\ref{lem:rank 0 center lemma}.
On the other hand, if $J=\{i,j\}$ is a non-edge pair, then 
 consider $g=g_ig_j$ with $g_i\in \ucompG_i$, $g_j\in \ucompG_j$.
For any $y\in \ucompG_i$, since $\ucompG_i$ is normal in $\ucompB_J$, 
 $y^d,y^g\in \ucompG_i$,  and since $\umap|_{\ucompG_i}$ is injective, we obtain $y^d=y^g=y^{g_i}$.
 As before we conclude that $g_i\in \ucompD_i$.
 Likewise $g_j\in \ucompD_j$ so that $g\in \ucompD$. 
The conclusion that $x\in Z(\ucompG)\cap \ucompD$ follows as above.  
\epf

\bco\label{cor:ker phi = central and diagonal}
For every $J\sbe I$, we have 
 $\ucompK_J\le Z(\ucompG)\cap \ucompD$.
In particular, $\ker\umap\le Z(\ucompG)\cap \ucompD$.
\eco
\bpf
This follows from an inductive application of Lemma~\ref{lem:tK=ker phi} combined with Lemma~\ref{lem:tK rank 2}. 
\epf

\ble\label{lem:phi is onto on central diagonal group}
$$Z(\ucompG)\cap\ucompD=\umap^{-1}(Z(\compG)\cap \compD).$$
\ele
\bpf
The inclusion `$\subseteq$' follows easily as $\umap(Z(\ucompG))\subseteq Z(\compG)$, since $\umap$ is onto,  and $\umap(\ucompD)\subseteq \ocompD$ by definition.
To see the converse inclusion, suppose that $d\in Z(\compG)\cap \ocompD$.
Then, since $Z(\compG)\cap \compD=Z(\compG)\cap \compD_{\theta}$ by assumption, there exists $\tilde{d}\in \ucompD$ with $\umap(\tilde{d})=d$.
By Lemma~\ref{lem:rank 0 center lemma} $\tilde{d}\in Z(\ucompG)\cap \ucompD$.
By Corollary~\ref{cor:ker phi = central and diagonal} $\umap^{-1}(d)=\tilde{d}\cdot \ker \umap \subseteq Z(\ucompG)\cap \ucompD$ as well.
\epf

\ble\label{lem:ZGtheta cap D in ZG}
$Z(\compG)\cap \ocompD= Z(\compL)\cap \compG$.
\ele
\bpf
We first claim that  $Z(\compG)\cap \compD\le Z(\compL)\cap \compG$.
Recall $\ocompD=\Stab_\compL(c,c^\theta)$, where $(c,c^\theta)\in \De^\theta$.
Let $g\in Z(\compG)\cap \ocompD$, then clearly $gc=c$ and, for any $h\in \compG$ we have 
 $g hc= hgc=hc$. By Corollary~\ref{cor:BJ flag transitive}, $g$ fixes all of $\De^\theta$.
 
We now claim that $g$ fixes all of $\De$.
Let $\pi$ be any panel on $c$. Then, 
 $\{d\in \pi\mid (d,d^\theta)\in \De^\theta\}=\pi-\{\proj^*(c^\theta)\}$.
It follows that $g$ fixes all chambers of $\pi$ and hence 
 $g$ fixes $E_1(c)\cup \{c^\theta\}$.
By Theorem 1 of~\cite[\S 4.2]{Ti1992}, $g$ acts as the identity on $\De$.
Therefore, by Lemma 3.4 of~\cite{Cap2007}, $g\in Z(\compL)$, proving the claim.

Finally, note that we have inclusions
$$ Z(\compL)\cap \compG=Z(\compL)\cap\ocompD\cap\compG\subseteq Z(\compG)\cap \ocompD\subseteq Z(\compL)\cap \compG.$$
The first equality follows from Theorem~\ref{thm:KM center}. 
As for the second equality follows from the inequality just obtained.
\epf

\bth\label{thm:kerphi = ZtGcaptD}
{Suppose that $\compL$ is the adjoint Kac-Moody group of type $\Gamma$, that is, $Z(\compL)=1$.}
Then, $\ker\umap=Z(\ucompG)\cap \ucompD$.
\eth
\bpf
The inclusion $\subseteq$ is Corollary~\ref{cor:ker phi = central and diagonal}.
To establish the converse inclusion, note that by Lemmas~\ref{lem:phi is onto on central diagonal group}~and~\ref{lem:ZGtheta cap D in ZG}, we have $\umap(Z(\ucompG)\cap \ucompD)\subseteq Z(\compG)\cap \ocompD=Z(\compL)\cap \ocompD=1$.
\epf

\medskip
\paragraph{Proof of Theorem~\ref{thm:concrete completions}.4}
The fact that $\ker\umap=Z(\ucompG)\cap\ucompD$ is Theorem~~\ref{thm:concrete completions}.
The equality $\ker\umap=d^\amG(\ker\comop)$ is Theorem~\ref{thm:comop characterizes center}.

\medskip
We finish this section with an observation showing that, although it isn't a group of Kac-Moody type,  the twisted group has the notion of a Weyl group.
Recall $W=\ocompN/\ocompD$, where $\ocompN=\Stab_{\compL}(\Sigma_\pm)$ and $\ocompD=\Stab_\compL((c,c^{\theta}))$.
Since $\Sigma_\pm=\Sigma(c,c^{\theta})$, $\theta$ acts naturally on $W$, and we let 
\begin{align*}
W^\theta&=\{w\in W \mid	w^\theta = w\},\\
S^\theta &=\{ss^\theta \mid	s \in S\}.
\end{align*}

\ble\label{lem:Wtheta is Coxeter}
\begin{enumerate}
\item The pair $(W^{\theta},S^{\theta})$ is a Coxeter system whose diagram is the diagram induced on the $\theta$-orbits in $\hI$.

\item Define $N^{\theta}=\Stab_{\compLth}(\Sigma_\pm)$.
Then, we also have $W^{\theta}=N^{\theta}\ocompD/\ocompD\cong N^{\theta}/(N^{\theta}\cap \ocompD)=N^{\theta}/\compD^{\theta}$.
\end{enumerate}
\ele
\bpf 
The orbits of $\theta$ on $\hI$ are pairs $(i,i^\theta)$ with $i\in \hI$ such that $s_i$ and $s_{i^\theta}$ commute.
Thus $S^\theta$ is the set of longest words corresponding to the orbits of $\theta$ on $\hI$.
Therefore, by the main result of~\cite{Mu1993}, $(W^\theta,S^\theta)$ is a Coxeter system.

Consider the map $N \to W$ sending $n\mapsto nD$.
For each $s\in S$ select $n_s\in N$ so that $n_sD$ represents $s$.
Then, $\langle n_s n_s^\theta\mid s\in S\rangle\le N^\theta$.
By the first part of this lemma, this subgroup already maps surjectively to $W^\theta$. Since $N^\theta$ is fixed by $\theta$ also the image is fixed by $\theta$.
Now we have $N^\theta/(N^\theta\cap \ocompD)\cong W^\theta$.
Note that $N^\theta\cap \ocompD=\compD_\theta$.
\epf

\section{Examples}\label{section:examples}
In this section we use the methods developed in the preceding sections to compute $Z(\ucompG)\cap\ucompD$ for some  Curtis-Tits groups $\ucompG$ and we verify condition (D) for another group.
\subsection{The generalized Cartan matrix}
In this subsection we consider the situation of Subsection~\ref{subsec:centers in CT groups}.
\bde
Let $E$ be the subring with identity $1_E$ of $\End(\fk^*)$ generated by 
 all $\{\rho_{j,i}\mid i,j\in I, i\sim j\}$ with the natural operations given by 
  \begin{align*}
  (\alpha \oplus \beta)(x) &= \alpha(x)\cdot \beta(x)\\
  (\alpha \otimes \beta)(x) &= (\alpha\circ \beta)(x)
    \end{align*}
Here by $0_E$ we denote the map $x\mapsto 1_\fk$ for all $x\in \fk^*$ and $1_E=\id_{\fk^*}$ and the $\rho_{i,j}$ are those defined in Definition~\ref{dfn:comop}.
Henceforth, we shall impose the restriction that $\langle \rho_{j,i}\mid i\sim j, i,j\in I\rangle\le \Aut(\fk^*)$ is finite and commutative, so that in particular, the ring $E$ is commutative. This happens for instance in the case where $\fk$ is finite.
The {\em generalized Cartan matrix} is the matrix $\comop=(\comop_{i,j})_{i,j\in I}\in M_n(E)\le \End((\fk^*)^I)$ given by 
\begin{align*}
\comop_{i,j}=\begin{cases}
\id^{-2}&\mbox{ if } i = j\\
\rho_{j,i} & \mbox{ if }j\sim i \\
0 & \mbox{ if } i\not \sim j.
\end{cases}
\end{align*}
Theorem~\ref{thm:comop characterizes center} now says that $\ker\comop$ consists of all elements in $(\fk^*)^I$ mapping to central elements in $\compG$ under the map $d$.
\ede

\bre
Note that $\id^{-2}$ is the additive inverse of $\id\oplus \id$.
In particular, if all $\rho_{i,j}$ are the identity map, then the matrix representing $\comop$ is the negative of the Cartan matrix $(A_{i,j})_{i,j\in I}$, making the connection to the characterization of the center due to Steinberg concrete.
\ere
\subsection{Computation of centers}
We will compute the centers of the Curtis-Tits groups $\compG$ with the diagrams shown below.
\paragraph{Example 1}
Here we consider a group of Lie type with diagram $E_6$. By uniqueness of the amalgam we can assume that all $\rho_{i,j}$ are trivial.
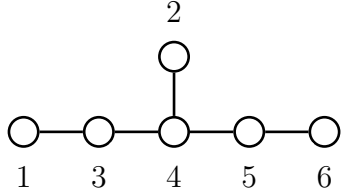
\begin{center}
  \begin{tikzpicture}[scale=.5]
  \tikzstyle{every node} = [draw, line width = 1pt, shape=circle]
    \node  [label=below:$1$] (one) at (0,0) {};
   \node  [label=above:$2$] (two) at  (4,2)  {};
   \node  [label=below:$3$] (three) at (2,0) {};
   \node  [label=below:$4$] (four) at (4,0) {};
   \node  [label=below:$5$] (five) at (6,0) {};
   \node  [label=below:$6$] (six) at (8,0)  {};
\foreach \from/\to in {one/three, three/four, four/five, five/six, two/four}
\path[draw, line width = 1pt]  (\from) -- (\to);
\path[draw, line width = 1pt, color=white] (three) .. controls (5,-3) .. (six);  \end{tikzpicture}
\end{center}
In the case of $E_6$, we find that $d=\prod_{i=1}^6 d_i(a_i)\in Z(\compG)$ if and only if the following conditions hold:
Setting
$a_1=a$, $a_6=b$, $a_2=c$, we get
$a_3=a^2$, $a_5=b^2$, $a_4=a^3=b^3=c^2$, $a_4^2=a^2b^2c$.
It then follows that $a^3=1=ab=c$, so the center is of order $3$ if and only if $\fk$ has a primitive third root of $1$ and trivial otherwise.

\paragraph{Example 2}
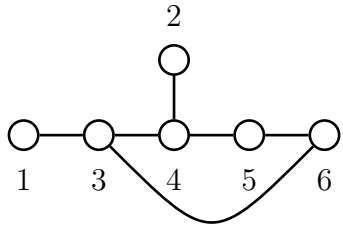
\begin{center}
 \begin{tikzpicture}[scale=.5]
  \tikzstyle{every node} = [draw, line width = 1pt, shape=circle]
    \node  [label=below:$1$] (one) at (0,0) {};
   \node  [label=above:$2$] (two) at  (4,2)  {};
   \node  [label=below:$3$] (three) at (2,0) {};
   \node  [label=below:$4$] (four) at (4,0) {};
   \node  [label=below:$5$] (five) at (6,0) {};
   \node  [label=below:$6$] (six) at (8,0)  {};
\foreach \from/\to in {one/three, three/four, four/five, five/six, two/four}
\path[draw, line width = 1pt]  (\from) -- (\to);
 \path[draw, line width = 1pt] (three) .. controls (5,-3) .. (six);
  \end{tikzpicture}
\end{center}

In this case we assume that  $\rho_{j,i}=1$ for all $1\le i\ne j\le 6$. This means that $\amG$ is orientable and $\ucompG$ is a split {group of Kac-Moody type}.
Assuming the field $\fk$ is finite, we can use a discrete logarithm to view the matrix $\comop$ as having coefficients in $\ZZ/(q-1)\ZZ$ and one computes that this matrix has determinant $-13$.
Therefore if $13$ is invertible modulo $q-1$, the center is trivial.
More precisely, we find that $d=\prod_{i=1}^6 d_i(a_i)\in Z(\compG)$ if and only if the following conditions hold:
Setting $a_1=a$, $a_2=c$, we get
$a_3=a^2$, $a_4=c^2$, $a_5=c^3a^{-2}$, $a_6=a^3c^{-2}$, $a_5^2=c^2a_6=a^3$ and similarly $a_6^2=c^3$.
It follows that $c^6=a^7$ and $a^6=c^7$ so that $ac=1=a^{13}$.
Therefore this $G$ has a non-trivial center of order $13$ if and only if $\fk$ contains a primitive $13$-th root of $1$.

\paragraph{Example 3}
Now suppose the amalgam has the same diagram as in Example 2, but  is given by 
$$\classmap_{j,i}=\begin{cases} \tau &\mbox{ if }(j,i)=(5,4)\\
 \id &\mbox{ otherwise }\end{cases}.$$
Then, $G$ is not orientable and $\rho_{5,4}=\rho_{4,5}\colon x\mapsto x^{-1}$. The matrix of the corresponding $\comop$ is that of the previous example replacing $\comop_{5,4}=\comop_{4,5}=-1$.
Now the determinant is $3$ and a computation similar to the previous one shows that $G$ now has a non-trivial center of order $3$ if and only if $\fk$ has a primitive $3$-rd root of $1$.
More precisely, $(a_1,\ldots,a_6)=(a,a^2,a^2,a,a^2,a^2)$, where $a^3=1$.

In case $\fk=\FF_q$ is finite of characteristic $p$, it may happen that $\classmap_{j,i}\colon x\mapsto x^{p^s}$ for some $s$.
In this case $\rho_{j,i}$ and $\rho_{i,j}$ are represented additively by $\comop_{j,i}=p^s$ and $\comop_{i,j}=-p^s$.
The determinant now should be computed modulo $q-1$, and in principle the center of $G$ can be computed as above.

\subsection{Verifying Condition~\texorpdfstring{\eqref{eqn:DGtheta=Dtheta}}{(D)}
}\label{example verifying condition D}
We consider the non-orientable amalgam $\amG$ of Example 3. The corresponding amalgam $\amL$ has the diagram shown below. The diagram is numbered using $\hI=\{i,i'\mid i\in \{1,2,3,4,5,6\}=I\}$, where $i'=\theta(i)$.
Note that $\amL$ corresponds to the map $\classmap\colon \hI\to \Aut(\fk)\times\langle\tau\rangle$ given by 
$$\classmap_{j,i}=\begin{cases} \tau &\mbox{ if }(j,i)\in \{(5,4'),(5',4)\}\\
 \id &\mbox{ otherwise }\end{cases}.
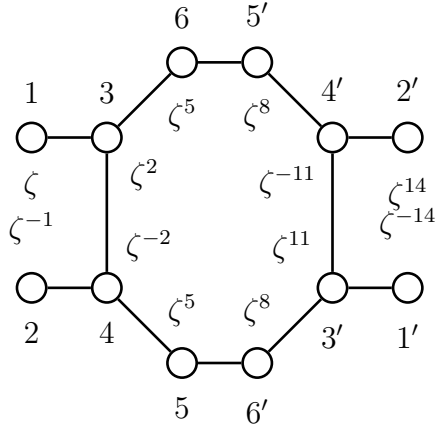
$$
\begin{center}
\begin{tikzpicture}[scale=.5]
  \tikzstyle{every node} = [draw, line width = 1pt, shape=circle]
    \node  [label=above:$1$, label=below:$\zeta$] (one) at (-5,2) {};
   \node  [label=below:$2$, label=above:$\zeta^{-1}$] (two) at  (-5,-2)  {};
   \node  [label=below:$4$, label=above right:$\zeta^{-2}$] (three) at (-3,-2) {};
   \node  [label=below:$5$, label=above:$\zeta^{5}$] (four) at (-1,-4) {};
   \node  [label=above:$6$, label=below:$\zeta^{5}$] (five) at (-1,4) {};
   \node  [label=above:$3$, label=below right:$\zeta^{2}$] (six) at (-3,2)  {};
    \node  [label=below:$1'$, label=above:$\zeta^{-14}$] (oneop) at (5,-2) {};
   \node  [label=above:$2'$, label=below:$\zeta^{14}$] (twoop) at  (5,2)  {};
   \node  [label=above:$4'$, label=below left:$\zeta^{-11}$] (threeop) at (3,2) {};
   \node  [label=above:$5'$, label=below:$\zeta^{8}$] (fourop) at (1,4) {};
   \node  [label=below:$6'$, label=above:$\zeta^{8}$] (fiveop) at (1,-4) {};
   \node  [label=below:$3'$, label=above left:$\zeta^{11}$] (sixop) at (3,-2)  {};
\foreach \from/\to in {one/six, two/three, three/six, three/four, five/six, four/fiveop,fourop/five, oneop/sixop,twoop/threeop,threeop/sixop,threeop/fourop,fiveop/sixop}
\path[draw, line width = 1pt]  (\from) -- (\to);
\end{tikzpicture}
\end{center}
A computation similar to the one in Example 3 reveals the following.
Any element  $a(\zeta)\in \kIK=\ker \comop$ is a sequence
 $(a_1,\ldots,a_6,a_{1'},\ldots,a_{6'})\in (\fk^*)^{\hI}$, where $a_i$ is the label of the node $i$ in the diagram and $\zeta^{39}=1$ in $\fk^*$.
 For example if we choose $a_1=\zeta$, then $a_{6'}=\zeta^{11}$.

By definition $a(\zeta)$ belongs to $\imnu\cap \kIK$ if and only if $\theta(a(\zeta))=a(\zeta)^{-1}$, that is
 $\zeta^{-1}=\zeta^{-14}$, which means that $\zeta^{13}=1$ in $\fk^*$.
 
Now, note that if $\mu_n$ is the set of $n$-th roots of $1$ in $\fk^*$, then the map $\mu_{39}\to \mu_{13}$ given by $x\mapsto x^{15}$ is surjective.
Noting that $\nu(a(\xi))=a(\xi^{15})$ for any $\xi$ a $39$-th root of $1$, we conclude that $(\imnu\cap \kIK)/\nu(\kIK)=0$, that is, condition~\eqref{eqn:DGtheta=Dtheta} is satisfied and $\compLth=\compG$ in this case.

\end{document}